\pdfoutput=1
\documentclass[11pt]{article}
\usepackage{color}
\usepackage{latexsym}
\usepackage{dsfont}
\usepackage{amssymb}
\usepackage{graphicx}
\usepackage{float}
\usepackage{amsmath, amsfonts,amssymb,theorem,euscript,array,enumerate,amsfonts,mathrsfs}
\usepackage[backend=bibtex, maxbibnames=10,maxcitenames=4, style=alphabetic]{biblatex}
\usepackage{caption}
\usepackage{subcaption}
\usepackage{nicefrac}
\usepackage{grffile}
\usepackage[noend]{algpseudocode}
\usepackage[ruled,vlined]{algorithm2e}
\usepackage[normalem]{ulem} %

\newtheorem{Remark}{Remark}[part]

\def \trans{^{\scriptscriptstyle{\intercal}}}

\newcommand{\cmt}[1]{\ignorespaces}
\usepackage{accents}
\newlength{\dhatheight}

\usepackage{enumitem}
\setlist[itemize]{leftmargin=*}

\newcommand*{\backin}{\rotatebox[origin=c]{-180}{$\in$}}

\def \Frac{\displaystyle\frac}

\def \b1{\bf{1}}

\def\argmax{\mathop{\rm argmax}}
\def\argmin{\mathop{\rm argmin}}

\newcommand{\red}[1]{\textcolor{red}{#1}}
\newcommand{\bl}[1]{\textcolor{blue}{#1}}

\def \A{\mathbb{A}}
\def \I{\mathbb{I}}
\def \N{\mathbb{N}}
\def \R{\mathbb{R}}
\def \M{\mathbb{M}}
\def \Z{\mathbb{Z}}
\def \E{\mathbb{E}}
\def \F{\mathbb{F}}
\def \G{\mathbb{G}}
\def \P{\mathbb{P}}
\def \Q{\mathbb{Q}}
\def \D{\mathbb{D}}
\def \T{\mathbb{T}}
\def \X{\mathbb{X}}

\def\esssup_#1{\underset{#1}{\mathrm{ess\,sup\, }}}

\def\argmin_#1{\underset{#1}{\mathrm{argmin\, }}}

\def \Ac{{\cal A}}
\def \Bc{{\cal B}}

\def \Fc{{\cal F}}

\def \Ic{{\cal I}}

\def \Nc{{\cal N}}

\def \Sc{{\cal S}}

\def \Wc{{\cal W}}

\def \eps{\varepsilon}

\def \ep{\hbox{ }\hfill$\Box$}

\def\Dt#1{\Frac{\partial #1}{\partial t}}

\def\reff#1{{\rm(\ref{#1})}}

\def\beqs{\begin{eqnarray*}}
\def\enqs{\end{eqnarray*}}
\def\beq{\begin{eqnarray}}
\def\enq{\end{eqnarray}}

\usepackage{hyperref}

\begin{document}

\title{
Deep neural networks algorithms for stochastic control problems 
on finite horizon: numerical applications\thanks{We are grateful to both referees for helpful comments and remarks.}
}

\author{
Achref \textsc{Bachouch}
\footnote{Department of Mathematics, University of Oslo, Norway. The  author's research is carried out with support of the Norwegian Research Council, within the research project Challenges in Stochastic Control, Information and Applications (STOCONINF), project number 250768/F20 \sf \href{mailto:achrefb at math.uio.no}{achrefb at math.uio.no}}
\quad
C\^ome \textsc{Hur\'e}
\footnote{LPSM, University Paris Diderot \sf \href{mailto:hure at lpsm.paris}{hure at lpsm.paris}}
\quad
Nicolas \textsc{Langren\'e}
\footnote{CSIRO Data61, RiskLab Australia \sf\href{mailto:Nicolas.Langrene at data61.csiro.au}{Nicolas.Langrene at data61.csiro.au}}
\quad
Huy\^en \textsc{Pham}
\footnote{LPSM, University  Paris-Diderot and CREST-ENSAE, {\sf  \href{mailto:pham at lspm.paris}{pham at lspm.paris}}
The work of this author is su\-pported by the ANR project CAESARS (ANR-15-CE05-0024), and also by FiME 
and the ``Finance and Sustainable Development'' EDF - CACIB Chair}
}

\maketitle 

\begin{abstract}
This paper presents several numerical applications of deep learning-based algorithms for discrete-time stochastic control problems in finite time horizon that have been introduced in  \cite{bacetal18_1}.  Numerical and comparative tests using {\sc TensorFlow}  
illustrate the performance of our different algorithms, namely control learning by performance iteration (algorithms NNcontPI and ClassifPI), control learning by hybrid iteration 
(algorithms Hybrid-Now and Hybrid-LaterQ),  
on the $100$-dimensional nonlinear PDEs examples from \cite{Ehanjen17} and on quadratic backward stochastic differential equations as in \cite{charic16}.   
We also performed tests on low-dimension control problems such as an option hedging problem in finance, as well as energy storage problems arising in the valuation of gas storage and in microgrid management.  Numerical results and compa\-risons to quantization-type algorithms Qknn, as an efficient algorithm to numerically solve low-dimensional control problems, are also provided.
\end{abstract}

\vspace{5mm}

\noindent \textbf{Keywords:} Deep learning, policy learning, performance iteration, value iteration, Monte\ Carlo, quantization.

\newpage

\section{Introduction} \label{secintro}

\setcounter{equation}{0} \setcounter{Assumption}{0}
\setcounter{Theorem}{0} \setcounter{Proposition}{0}
\setcounter{Corollary}{0} \setcounter{Lemma}{0}
\setcounter{Definition}{0} \setcounter{Remark}{0}

This paper is devoted to the numerical resolution of discrete-time stochastic control problem over a finite horizon. 
 The dynamics of the controlled state process $X$ $=$ $(X_n)_{n}$ valued in $\R^d$  is given by
 \beq \label{dynX}
 X_{n+1} &=& F(X_n,\alpha_n,\eps_{n+1}), \;\;\; n=0,\ldots,N-1, \; X_0 = x_0 \in \R^d, 
 \enq
 where $(\eps_n)_{n}$  is a sequence of i.i.d. random variables valued in some Borel space $(E,\Bc(E))$, and defined on some probability space $(\Omega,\Fc,\P)$ equipped with the filtration $\F$ $=$ $(\Fc_n)_n$ generated by the noise $(\eps_n)_n$ ($\Fc_0$ is the trivial $\sigma$-algebra), the control $\alpha$ $=$ $(\alpha_n)_{n}$ is an 
 $\F$-adapted process valued in  $\A$ $\subset$ $\R^q$, and $F$ is a measurable function from $\R^d\times\R^q\times E$ into $\R^d$ which is known by the agent. 
Given a running cost function $f$ defined on $\R^d\times\R^q$ and a terminal cost function $g$ defined on $\R^d$,  the cost functional associated with a control process  $\alpha$ is
 \beq \label{costJ}
 J(\alpha) &=& \E \left[ \sum_{n=0}^{N-1}  f(X_n,\alpha_n) +  g(X_N) \right]. 
 \enq 
 In this framework, we assume $f$ and $g$ to be known by the agent.
 The set $\Ac$ of admissible controls is the set of control processes $\alpha$ satisfying some integrability conditions  ensuring  that the cost functional $J(\alpha)$ is well-defined and finite. The control problem, also called Markov decision process (MDP), is formulated as  
\beq \label{defcontrol}
V_0(x_0) & := & \inf_{\alpha\in\Ac} J(\alpha),
\enq
and the goal is to find  an optimal control $\alpha^*$ $\in$ $\Ac$, i.e.,  attaining the optimal value: $V_0(x_0)$ $=$ $J(\alpha^*)$. 
Notice that problem  \reff{dynX}-\reff{defcontrol} may also be viewed as the time discretization of a continuous time stochastic control problem, in which case, $F$ is typically the Euler scheme for a controlled diffusion process.

It is well-known that the global dynamic optimization problem \reff{defcontrol} can be  reduced to local optimization problems via the dynamic programming (DP) approach, which allows to determine the value function in a backward recursion  by 
\beq
V_N(x) &=& g(x), \;\;\; x \in \R^d,  \nonumber \\
V_n(x) &=& \inf_{a\in\A} Q_n(x,a),   \label{DP}  \\
 \mbox{ with }  \; Q_n(x,a) & = &  f(x,a) + \E\big[ V_{n+1}(X_{n+1}) \big| X_n = x, \alpha_n = a  \big], \;\;\; (x,a) \in \R^d\times\A. \nonumber 
\enq
Moreover, when the infimum is attained in the DP formula \reff{DP} at any time $n$ by $a_n^*(x)$ $\in$ ${\rm arg}\min_{a \in \A} Q_n(x,a)$, we get an optimal control in feedback form (policy)  
given by:  $\alpha^*$ $=$  $(a_n^*(X_n^*))_n$ where $X^*$  is the Markov process defined by
\beqs
X_{n+1}^* &=& F(X_n^*,a_n^*(X_n^*),\eps_{n+1}), \;\;\;  n=0,\ldots,N-1, \;\; X_0^* = x_0. 
\enqs

The practical implementation of the DP formula may suffer from the curse of dimensionality and large complexity when the state space dimension $d$ and the control space dimension are high. In \cite{bacetal18_1}, we proposed algorithms relying on deep neural networks for approxima\-ting/learning the optimal policy and then eventually the value function by performance/policy iteration or hybrid iteration with Monte Carlo regressions now or later. 
This research led to three algorithms, namely algorithms NNcontPI,  Hybrid-Now  and Hybrid-LaterQ that are recalled in Section \ref{secalgo}, and which can be seen as a natural extension of actor-critic methods, developed in the reinforcement learning community for stationary stochastic problem (\cite{sutbar98}), to finite-horizon control problems. Note that for stationary control problem, it is usual to use techniques such as temporal difference learning, which relies on the fact that the value function and the optimal control do not depend on time, to improve the learning of the latter. Such techniques do not apply to finite horizon control problems.
In Section \ref{secnum}, we perform some numerical and comparative tests to  illustrate the efficiency of our different algorithms,  
on $100$-dimensional nonlinear PDEs examples as in \cite{Ehanjen17} and quadratic Backward Stochastic Differential equations as in \cite{charic16}, as well as on high-dimensional linear quadratic stochastic control problems.   
We present numerical results for an option hedging problem in finance,  and energy storage problems arising  in the valuation of gas storage and in microgrid management. 
Numerical results and comparisons to quantization-type algorithms Qknn, introduced in this paper as an efficient algorithm to numerically solve low-dimensional control problems, are also provided.
Finally, we conclude in Section \ref{secconclu} with some comments about possible  extensions and improvements of our algorithms.

\section{Algorithms} \label{secalgo}

\setcounter{equation}{0} \setcounter{Assumption}{0}
\setcounter{Theorem}{0} \setcounter{Proposition}{0}
\setcounter{Corollary}{0} \setcounter{Lemma}{0}
\setcounter{Definition}{0} \setcounter{Remark}{0}

We introduce in this section four  neural network-based algorithms for  solving  the discrete-time stochastic control problem \reff{dynX}-\reff{defcontrol}. The convergence of these algorithms have been analyzed in detail in our companion paper \cite{bacetal18_1}, and for self-contained purpose, we recall in this section the description of these algorithms and the convergence results. 
We also introduce at the end of this section a quantization and $k$-nearest-neighbor-based algorithm (Qknn) that will be used as benchmark when testing our algorithms on low-dimensional control problems.

We are given a class of deep neural networks (DNN)  for the control policy represented by the parametric functions $x$ $\in$ $\R^d$ $\mapsto$ $A(x;\beta)$ $\in$ $\A$, with parameters $\beta$ $\in$ $\R^q$,  and a class of DNN  for the value function represented by the parametric functions: $x$ $\in$ 
$\R^d$ $\mapsto$ $\Phi(x;\theta) \in \R$, with parameters $\theta$ $\in$ $\R^p$. Recall that these DNN functions $A$ and $\Phi$ 
are compositions of linear combinations and nonlinear activation functions, see \cite{good16}.

Additionally, we shall be given a sequence of probability measures on the state space $\R^d$, that we call training measure and denoted $(\mu_n)_{n=0}^{N-1}$, which should be seen as dataset providers to learn the optimal strategies and the value functions at time $n=0,\ldots,N-1$. 

\begin{Remark}[Training sets design]
	\label{rk:trainingMeasure}
{\rm The choice of the training sets is critical for numerical efficiency. This problem has been largely investigated in the reinforcement learning community, notably with multi-armed bandits algorithms \cite{Auer2002}, and more recently in the numerical probability literature, see \cite{Ludkovski2019}, but remains a  challenging issue.  
Here, two cases are considered for  the choice of the training measure $\mu_n$ used to generate the training sets on which the estimates at time $n$ will be computed. 
The first one is a knowledge-based selection, relevant when the controller knows with a certain degree of confidence where the process has to be driven in order to optimize her cost functional. The second case is when the controller has no idea where or how to drive the process to optimize the cost functional. 

\vspace{2mm}

\noindent {\sc (1) Exploitation only strategy} 

\vspace{1mm}

\noindent In the knowledge-based setting, there is no need for exhaustive and expensive (in time mainly) exploration of the state space, and the controller can take a training measure $\mu_n$ that assigns more points in the region of the state space that is likely to be visited by the optimally-driven process. 

In practice, at time $n$, assuming we know that the optimal process is likely to lie in a region $\mathcal{D}$, we choose a training measure in which the density assigns a lot of weight to the points of $\mathcal{D}$, for example $\mathcal{U}(\mathcal{D})$, the uniform distribution in $\mathcal{D}$.

\vspace{2mm}

\noindent {\sc (2) Explore first, exploit later}
When the controller has no idea where or how to drive the process to optimize the cost functional, we suggest to build the training measures as empirical measures of the process, driven by  estimates of the optimal control computed using alternative methods.
\begin{itemize}
	\item[(i)] {\it Explore first:} Use an alternative method to obtain good estimates of the optimal strategy. In high-dimension: one can for example think of approximating the control at all time by neural network, and obtain a good estimate of the optimal control by performing a global optimization of the function:
	\[
	J(\theta_0,\ldots,\theta_{N+1}) := \E\left[ \sum_{n=0}^{N-1} f(X_n,A(X_n;\theta_n)) + g(X_N)   \right],
	\]
	where $X$ is the process controlled by the feedback control $A(.;\theta_n)$ at time $n$.\\
	\item[(ii)] {\it Exploit later:}
	Take the training measures $\mu_n:=\P_{X_n}$, for $ n=0,\ldots,N-1$, where $X$ is driven using the optimal control estimated in step (i); and apply the procedure (1).
Such an idea has been recently exploited in \cite{kouetal18}.  
\end{itemize}
}
\end{Remark}	
	
\begin{Remark}[Choice of Neural Networks]
{\rm
	Unless otherwise specified, we use feed-forward Neural Networks with two or three hidden layers and $d$+10 neurons per hidden layer, since we noticed empirically that these parameters were enough to approximate the relatively smooth objective functions considered here. We tried sigmoid, tanh, ReLU and ELU activation functions and noticed that ELU is most often the one providing the best results in our applications. We normalize the input data of each neural network in order to speed up the training of the latter.
 }
\ep
\end{Remark}

\begin{Remark}[Neural Networks Training]
{\rm
We use the Adam optimizer, as implemented in {\sc TensorFlow}, with initial learning-rate set to 0.001 or 0.005, which are the default va\-lues in {\sc TensorFlow}, to train by gradient-descent the optimal strategy and the value function defined in the algorithms described later. {\sc TensorFlow} takes care of the Adam gradient-descent procedure by automatic differentiation when the function to optimize is an expectation of {\sc TensorFlow} functions, such as the usual differentiable activation functions $\sin, \log, \exp$ but also popular non-differentiable activation functions such as ReLu: $x \mapsto \max(0,x)$.

 In order to force the weights and  biases of the neurons to stay small, we use an $\mathbb{L}^2$ regularization with parameter mainly set to 0.01, but the value can change in order to make sure that the regularization term is neither too strong or too weak when added to the loss when training neural networks.
	
We consider a large enough number of mini-batches of size 64 or 128 for the training, depending essentially empirically on the dimension of the problem. We use at least 10 epochs\footnote{We denote by epoch one pass of the full training set.} and stop the training when the loss computed on a validation set of size 100 stops decreasing. We noticed that taking more than one epoch really improves the quality of the estimates.
}
\ep
\end{Remark}

\begin{Remark}[Constraints] \label{remcons}
{\rm The proposed algorithms can deal with state and control constraints at any time, which is useful in several applications: 
\beqs
(X_n^\alpha,\alpha_n) & \in & \Sc \;\;\; a.s., \;\; n  \in \N, 
\enqs
where $\Sc$ is some given subset of  $\R^d\times\R^q$. In this case, in order to ensure that the set of admissible controls is not empty, we assume that the sets
\beqs
\A(x) & := & \Big\{ a \in \R^q: (F(x,a,\eps_{1}),a)  \in \Sc  \; a.s. \Big\} 
\enqs
are non empty for all $x$ $\in$ $\Sc$, and the DP formula now reads
\beqs
V_n(x) &=& \inf_{a\in \A(x)} \big[ f(x,a) +  P^a V_{n+1}(x)\big], \;\;\; x \in \Sc.  
\enqs
From a computational point of view, it may be more convenient to work with unconstrained state/control variables, hence by relaxing the state/control constraint and introducing into the running cost a penalty function $L(x,a)$: $f(x,a)$ $\leftarrow$ $f(x,a) + L(x,a)$, and $g(x)$ $\leftarrow$ $g(x) + L(x,a)$.  
For example, if the constraint set $\Sc$ is in the form: 
$\Sc$ $=$ $\{ (x,a) \in \R^d\times\R^q: h_k(x,a) =  0, k=1,\ldots,p, \; h_k(x,a) \geq 0, k=p+1,\ldots,q\}$, for some functions $h_k$, then one can take as  penalty functions:
\beqs
L(x,a) &=& \sum_{k=1}^p \mu_k |h_k(x,a)|^2 +  \sum_{k=p+1}^q \mu_k \max(0, - h_k(x,a) ).  
\enqs 
where $\mu_k$ $>$ $0$ are penalization  coefficients (large in practice). 
}
\ep
\end{Remark}

\subsection{Control Learning by Performance Iteration}

We present in this section Algorithm \ref{algo:nncontPI}, which combines an optimal policy estimation by neural networks and the dynamic programming principle. We rely on the performance iteration procedure, i.e. paths are always recomputed up to the terminal time $N$.

\subsubsection{Algorithm NNContPI}

Our first algorithm, referred to as NNContPI,  
is well-designed for control problems with continuous control space such as $\R^q$ or a ball in $\R^q$. The main idea is:
\begin{enumerate}
	\item Represent the controls at time $n=0,\ldots, N-1$ by neural networks in which the activation function for the output layers takes values in the control space. For example, one can take the identity function as activation function for the output layer if the control space is $\R^q$; or the sigmoïd function if the control space is $[0,1]$.
	\item Learn sequentially in  time, and in a backward way, the optimal parameters $\hat{\beta}_n$ for the representation of the optimal control. In particular, notice that the learning of the optimal control at time $n$ highly relies on the accuracy of the estimates of the optimal controls at time $k=n+1,\ldots,N-1$, computed previously.
\end{enumerate}

\begin{algorithm}
 	\caption{NNContPI}
 	\label{algo:nncontPI}
 	\textbf{Input:} the training distributions $(\mu_n)_{n=0}^{N-1}$\;
 	\textbf{Output:} estimates of the optimal strategy $(\hat a_n)_{n=0}^{N-1}$\;
 	\For{$n$ $=$ $N-1,\ldots,0$}{
 		Compute
 \beq 
 \hat \beta_n & \in  &    \argmin_{\beta \in \R^q} \E \left[ f\big(X_n^{},A(X_n^{};\beta)\big)+  \sum_{k=n+1}^{N-1} 
 f\big(  X^{\beta}_k, \hat{a}_k \big( X^{\beta}_k \big)\big)  + g\big(  X^{\beta}_N \big)\right]\; \nonumber\\
 \label{computa}
 \enq	
 where $X_n \sim \mu_n$ and where $\big(X^{\beta}_k\big)_{k=n+1}^N$ is defined by induction as:%
 \begin{equation*}
 \left\{
 \begin{array}{ccl}
 X^{\beta}_{n+1} &=& F\big( X^{}_n, A\big(X^{}_n;\beta \big), \eps_{n+1}^{} \big) \\
 X^{\beta}_{k+1} &=& F\big( X^{\beta}_{k}, \hat a_k\big(X^{\beta}_{k}\big), \eps_{k+1}^{} \big), \;\;\; \mbox{  for } k= n+1, \ldots, N-1.
 \end{array}
 \right.
 \end{equation*}
Set $\hat a_n  =  A(.;\hat\beta_n)$. 

\Comment{$\hat a_n$ is the estimate of the optimal policy at time $n$}
 }
 \end{algorithm}

\subsubsection{Algorithm ClassifPI}

In the special case where the control space $\A$ is finite, i.e., Card$(\A)$ $=$ $L$ $<$ $\infty$ with  $\A$ $=$ $\{a_1,\ldots,a_L\}$, a classification method can be used: 
 consider a DNN that takes state $x$ as input and returns a probability vector $p(x;\beta)$ $=$ $(p_\ell(x;\beta))_{\ell=1}^L$  with parameters $\beta$. Such a usual DNN can be build using $k$ hidden layers with ReLu activation functions, an output layer with $L$ neurons, and a Softmax\footnote{The Softmax function is defined as follows: $x \mapsto \left( \frac{e^{\beta_1 x}}{\sum_{k=1}^{L} e^{\beta_k x}}, \ldots, \frac{e^{\beta_1 x}}{\sum_{k=1}^{L} e^{\beta_k x}} \right)$ where $\beta_1,\ldots,\beta_L$ are part of the parameters that will be learned by gradient-descent.} activation function for the output layer. 
 Algorithm \ref{algo:classifPI}, presented below, is based on this idea, and is called ClassifPI. 
 
\begin{algorithm}
 	\caption{ClassifPI}
 	\label{algo:classifPI}
 	\textbf{Input:} the training distributions $(\mu_n)_{n=0}^{N-1}$\;
 	\textbf{Output:} estimates of optimal strategies $(\hat a_n)_{n=0}^{N-1}$ and probabilities $p_l(.;\hat{\beta}_n)$\;
 	\For{$n$ $=$ $N-1,\ldots,0$}{
 		Represent the discrete control at time $n$ by neural network with parameter $\beta_n$:
 		\beq
 		a_n(x) &=&  a_{\ell_n(x)} \;  \mbox{ with } \;  \ell_n(x) \; \in \; \argmax_{\ell = 1,\ldots,L} p_\ell(x;\beta_n), \nonumber
 		\enq 
 		and compute the optimal parameter:
 		\beq
 		\hat\beta_n & \in & \argmin_{\beta\in\R^q}  \E \left[  \sum_{\ell=1}^L p_\ell(X_n^{};\beta) \Big( f(X_n^{},a_\ell)   
 		+  \sum_{k=n+1}^{N-1} f\big(X_k^{\ell},\hat a_k(X_k^{\ell})\big) \; + \; g(X_N^{\ell})  \Big) \right], \nonumber \\
 		\label{computPI}
 		\enq
 		where $X_n\sim \mu_n$ on $\R^d$,  $X_{n+1}^{\ell}$ $=$ $F(X_n^{},a_\ell,\eps_{n+1}^{})$, $X_{k+1}^{\ell}$ $=$ $F(X_k^{\ell},\hat a_k(X_k^{\ell}),\eps_{k+1}^{})$, 
 		for $k=n+1,\ldots,N-1$ and $\ell$ $=$ $1,\ldots,L$\;
 		Set $ \hat a_n(.) = a_{\hat \ell_n(.)} \;  \mbox{ with } \;  \hat\ell_n(x) \; \in \; \argmax_{\ell = 1,\ldots,L} p_\ell(x;\hat\beta_n)$\;
		
		\Comment{$\hat a_n$ is the estimate of the optimal policy at time $n$}
 	}
 \end{algorithm}

Note that, when using Algorithms \ref{algo:nncontPI} and \ref{algo:classifPI}, the estimate of the optimal strategy at time $n$ highly relies on the estimates of the optimal strategy at time $n+1,\ldots,N-1$, that have been computed previously. In particular, the practitioner who wants to use Algorithms \ref{algo:nncontPI} and \ref{algo:classifPI} needs to keep track of the estimates of the optimal strategy at time $n+1,\ldots,N-1$ in order to compute the estimate of the optimal strategy at time $n$.

\begin{Remark} \label{rk:approxExpectations}
{\rm  	In practice, for $n=N-1,...,0$, one should minimize the expectations \eqref{computa} and \eqref{computPI} by stochastic gradient-descent, where mini-batches of finite number of paths 
$(X_k^\beta)_{k=n+1}^N$ are generated by drawing independent samples under $\mu_n$ for the initial position at time $n$, and independent samples under $\varepsilon_{k}$, for $k=n+1,\ldots,N$.
The convergence of Algorithms 1 and 2 is analyzed in \cite{bacetal18_1} in terms of the error approximation of the optimal control by neural networks, and in terms of the estimation error 
 by stochastic gradient descent methods, see their Theorem 4.7.  
}
\ep
 \end{Remark}

\subsection{Control and value function learning by double DNN}   

We present in this section two algorithms, which in contrast with Algorithms \ref{algo:nncontPI} or \ref{algo:classifPI},  only  keep track of the estimates of the value function and optimal control at time $n+1$ in order to build an estimate of the value function and optimal control at time $n$.

\subsubsection{Regress Now (Hybrid-Now)}

The  Algorithm \ref{algo:Hybrid-Now}, refereed to as Hybrid-Now, combines optimal policy estimation by neural networks and dynamic programming principle, and relies on  an hybrid procedure between value and performance iteration.

 \begin{algorithm}
 	\caption{Hybrid-Now}
 	\label{algo:Hybrid-Now}
 	\textbf{Input:} the training distributions $(\mu_n)_{n=0}^{N-1}$\;
 	\textbf{Output:} \\
 	-- estimate of the optimal strategy $(\hat a_n)_{n=0}^{N-1}$\;
 	-- estimate of the value function $(\hat V_n)_{n=0}^{N-1}$\;
 Set $\hat V_N$ $=$ $g$\;
 	\For{$n$ $=$ $N-1,\ldots,0$}{
 Compute: 
\beq
\hat \beta_n  & \in  &    \argmin_{\beta \in \R^q} \E \Big[ f\big(X_n^{},A(X_n^{};\beta)\big)+ \hat{V}_{n+1}(X_{n+1}^{\beta}) \Big] \label{computanow}
\enq
where $X_n \sim \mu$, and 
$X_{n+1}^{\beta}$ $=$ $F\big(X_n^{},A(X_n^{};\beta)\big),\eps_{n+1}^{})$\;
 		Set $ \hat a_n  =  A(.;\hat\beta_n) $; \Comment{$\hat a_n$ is the estimate of the optimal policy at time $n$}\\
 Compute
 \beq
 \hat\theta_n & \in & \argmin_{\theta\in\R^p} \E  \left[ \Big(( f(X_n^{},\hat a_n(X_n^{}))+ \hat{V}_{n+1}(X_{n+1}^{\hat\beta_n}) - \Phi(X_n^{};\theta) \Big)^2 \right].  \label{computv}
 \enq 
 Set $\hat V_n = \Phi(.;\hat\theta_n) $; \Comment{$\hat V_n$ is the estimate of the value function at time $n$}
 	}
 \end{algorithm}

\begin{Remark}
{\rm	One can combine different features from Algorithms \ref{algo:nncontPI}, \ref{algo:classifPI} and \ref{algo:Hybrid-Now} to solve specific problems, as it has been done for example in Section \ref{sec:sgm}, where we designed Algorithm \ref{algo:ClassifHybrid} to solve a smart grid management problem.
}
\ep
\end{Remark}

\subsubsection{Regress Later and Quantization (Hybrid-LaterQ)}   
\label{sec:Hybrid-LaterQ}

The   Algorithm \ref{algo:Hybrid-LaterQ}, called Hybrid-LaterQ,  combines regress-later and quantization methods to build estimates of the value function. The main idea behind Algorithm \ref{algo:Hybrid-LaterQ} is to first interpolate the value function at time $n+1$ by a set of basis functions, which is in the spirit of the regress-later-based algorithms, and secondly regress the interpolation at time $n$ using quantization. The usual regress-later approach requires the ability to compute closed-form conditional expectations, which limits the stochastic dynamics and regression bases that can be considered. The use of quantization avoids this limitation and makes the regress-later algorithm more generally applicable.

Let us first recall the basic ingredients of quantization. We denote by $\hat\eps$ a $K$-quantizer of  the $\R^d$-valued random variable $\eps_{n+1}$ $\sim$ $\eps_1$ (typically a Gaussian random variable), 
	that is a discrete random variable on a grid $\Gamma$ $=$ $\{e_1,\ldots,e_K\}$ $\subset$ $(\R^d)^K$ defined by 
	\beqs
	\hat\eps &=& {\rm Proj}_\Gamma(\eps_1) \; := \; \sum_{\ell=1}^K e_\ell 1_{\eps_1 \in C_\ell(\Gamma)},
	\enqs
	where $C_1(\Gamma)$, $\ldots$, $C_K(\Gamma)$ are Voronoi tesselations of $\Gamma$, i.e., Borel partitions of the Euclidian space $(\R^d,|.|)$ satisfying
	\beqs
	C_\ell (\Gamma) & \subset & \Big\{ e \in \R^d: |e-e_\ell | \; = \; \min_{j =1,\ldots,K} |e- e_j  | \Big\}. 
	\enqs
	The discrete law of $\hat\eps$ is then characterized by 
	\beqs
	\hat p_\ell & :=& \P[ \hat\eps = e_\ell ] \; = \; \P[ \eps_1 \in C_\ell (\Gamma) ], \;\;\; \ell=1,\ldots,K.  
	\enqs
	The grid points $(e_\ell )$ which minimize the $L^2$-quantization error $\| \eps_1 - \hat\eps\|_{_2}$ lead to the so-called optimal $K$-quantizer, and can be obtained by a stochastic gradient descent method, known as Kohonen algorithm or competitive learning vector quantization (CLVQ) algorithm, which also provides  as a byproduct an estimation of the associated weights $(\hat p_\ell )$. We refer to \cite{pagetal04} for a description of the algorithm, and mention that for the normal distribution, the optimal grids and the weights of the Voronoi tesselations are precomputed on the website 
	\url{http://www.quantize.maths-fi.com}.

\vspace{3mm}

\begin{algorithm}[H]
	\caption{Hybrid-LaterQ}
	\label{algo:Hybrid-LaterQ}
	\textbf{Input:} \\
	-- the training distributions $(\mu_n)_{n=0}^{N-1}$\;
	-- The grid $\{e_1,\ldots,e_K\}$ of $K$ points in $\R^d$, with weights $p_1,\ldots,p_K$ for the quantization of the noise $\eps_n$\;
	\textbf{Output:} \\
	-- estimate of the optimal strategy $(\hat a_n)_{n=0}^{N-1}$\;
	-- estimate of the value function $(\hat V_n)_{n=0}^{N-1}$\;
	Set $\hat V_N$ $=$ $g$\;
	\For{$n$ $=$ $N-1,\ldots,0$}{
		Compute: 
		\beq 
		\hat \beta_n  & \in  &    \argmin_{\beta \in \R^q} \E \Big[ f\big(X_n^{},A(X_n^{};\beta)\big)+ \hat{V}_{n+1}(X_{n+1}^{\beta}) \Big] \label{computalater}
		\enq
		where $X_n \sim \mu_n$, and 
		$X_{n+1}^{\beta}$ $=$ $F\big(X_n^{},A(X_n^{};\beta)\big),\eps_{n+1}^{})$\;
		Set $ \hat a_n  =  A(.;\hat\beta_n) $; \Comment{$\hat a_n$ is the estimate of the optimal policy at time $n$}\\
		Compute
		\beq
		\hat\theta_{n+1} & \in & \argmin_{\theta\in \mathbb{R}^p} \E  \left[  \left( \hat{V}_{n+1}(X_{n+1}^{\hat\beta_n}) - \Phi(X_{n+1};\theta) \right)^2 \right]  \label{computvlater}
		\enq 
		and set $\tilde V_{n+1} = \Phi(.;\hat\theta_{n+1})$\; \Comment{interpolation at time $n+1$}\\
		Set 
\beqs
\hat V_n(x) & =&  f(x,\hat{a}_n(x))+ \sum_{\ell=1}^K p_\ell \tilde{V}_{n+1}\big( F(x,\hat{a}_n(x),e_\ell )\big);
\enqs
\Comment{$\hat V_n$ is the estimate by quantization of the value function at time $n$}\\
	}
\end{algorithm}

\vspace{3mm}

Quantization is mainly used in Algorithm \ref{algo:Hybrid-LaterQ} to efficiently approximate the expectations: recalling the dynamics \reff{dynX}, the conditional expectation operator for any functional $W$ is equal to
	\beqs
	P^{\hat a_n^M(x)} W(x) \; = \; \E\big[ W(X_{n+1}^{\hat a_n^M}) | X_n = x \big] &=& \E \big[ W(F(x,\hat a_n^M(x),\eps_1)) \big], \;\; x \in \R^d, 
	\enqs
	that we shall approximate analytically by quantization via: 
	\beqs
	\widehat P^{\hat a_n^M(x)} W(x) &:=& \E \big[ W(F(x,\hat a_n^M(x),\hat\eps)) \big] \; = \; \sum_{\ell=1}^K \hat p_\ell  W \left(F(x,\hat a_n^M(x),e_\ell ) \right). 
	\enqs

Observe that the solution to \eqref{computvlater} actually provides a neural network $\Phi(.;\hat\theta_{n+1})$ that interpolates $\hat{V}_{n+1}$. Hence the Algorithm \ref{algo:Hybrid-LaterQ} contains an interpolation step, and moreover, any kind of distance in $\R^d$ can be chosen as a loss to compute $ \hat \theta_{n+1}$. In \reff{computvlater}, we decide to take the $\mathbb{L}^2$-loss, mainly because it is the one that worked the best in our applications.

\begin{Remark}[Quantization]
{\rm 
	In dimension 1, we used the optimal grids and weights with $K=21$ points, to quantize the reduced and centered normal law $\mathcal{N}(0,1)$; and took 100 points to quantize the reduced and centered normal law in dimension 2, i.e. $\mathcal{N}_2(0,1)$. All the grids and weights for the optimal quantization of the normal law in dimension $d$  are  available in \url{http://www.quantize.maths-fi.com} for $d=1,\ldots,100$.
}
\ep 	
\end{Remark}

\subsubsection{Some remarks on Algorithms \ref{algo:Hybrid-Now}  and \ref{algo:Hybrid-LaterQ}}
\label{sec:pretrain}

As in Remark \ref{rk:approxExpectations}, all the expectations written in our pseudo-codes  in Algorithm \ref{algo:Hybrid-Now}  and \ref{algo:Hybrid-LaterQ} should  be approximated by empirical mean using a finite training set. The convergence of these algorithms has been analyzed in \cite{bacetal18_1} in terms of the  approximation error 
of the optimal control  and value function by  neural networks, in terms of the estimation error  by stochastic gradient descent methods, and in terms of the quantization error (for Algorithm \ref{algo:Hybrid-LaterQ}, see their Theorems 4.14 and 4.19).

Algorithms \ref{algo:Hybrid-Now} or \ref{algo:Hybrid-LaterQ} are quite efficient to use in the usual case where the value function and the optimal control at time $n$ are very close to the value function and the optimal control at time $n+1$, which happens e.g. when the value function and the optimal control are approximations of the time discretization of a continuous in time value function and an optimal control. In this case, it is recommended to follow this two-step procedure:
	 \begin{itemize}
	 	\item[(i)] initialize the parameters (i.e. weights and bias) of the neural network approximations of the value function and the optimal control at time $n$ to the ones of the neural network approximations of the value function and the optimal control at time $n+1$.
	 	\item[(ii)] take a very small learning rate parameter, for the Adam optimizer, that guarantees the stability of the parameters' updates from the gradient-descent based learning procedure.
	 \end{itemize} 
Doing so, one obtains stable estimates of the value function and optimal control, which is desirable. We highlight the fact that this stability procedure is applicable here since the stochastic gradient descent method  benefits from good initial guesses of the parameters to be optimized. It is an advantage compared to alternative methods proposed in  the literature, such as classical polynomial regressions. 

\subsection{Quantization with k-nearest-neighbors (Qknn-algorithm)}   

Algorithm \ref{algo:Qknn} presents  the pseudo-code of an algorithm based on the quantization and $k$-nearest neighbors methods,  called Qknn, which will be the benchmark in all the low-dimensional control problems that will be considered in Section \ref{secnum} to test NNContPI, ClassifPI, Hybrid-Now and Hybrid-Later. Also, comparisons of Algorithm \ref{algo:Qknn} to other well-known algorithms on various control problems in low-dimension are performed in \cite{BalHurPha18}, which show in particular that Algorithm \ref{algo:Qknn} works very well to solve low-dimensional control problems. Actually, in our experiments, Algorithm \ref{algo:Qknn} always outperforms the other algorithms based either on regress-now or regress-later methods whenever the dimension of the problem is low enough for Algorithm  \ref{algo:Qknn} to be feasible.

\vspace{2mm}

As done in Section \ref{sec:Hybrid-LaterQ}, we  consider a $K$-optimal quantizer of the noise $\eps_n$, i.e. a discrete random variable $\hat\eps_n$ valued in a  grid $\{e_1,\ldots,e_K\}$ of $K$ points in $E$, and with weights $p_1,\ldots,p_K$. We also consider grids $\Gamma_n$, $n=0,\ldots,N$ of points in $\R^d$, which are assumed to properly cover the region of $\R^d$ that is likely to be visited by the optimally driven process $X$ at time $n=0,\ldots,N-1$. These grids can be viewed as  samples of  well-chosen training distributions where more points are taken in the region that is likely to be visited by the optimally driven controlled process (see Remark \ref{rk:trainingMeasure} for details on the choice of the training measure).

\vspace{3mm}

\begin{algorithm}[H]
	\caption{Qknn}
	\label{algo:Qknn}
	\textbf{Input:}\\
	-- Grids $\Gamma_k$, $k=0,\ldots,N$ in $\R^d$\;
	--  Grid $\{e_1,\ldots,e_K\}$ of  $K$ points in $E$, with weights $p_1,\ldots,p_K$ for the quantization of $\eps_n$ \\
	\textbf{Output:} \\
	-- estimate of the optimal strategy $(\hat a_n)_{n=0}^{N-1}$\;
	-- estimate of the value function $(\hat V_n)_{n=0}^{N-1}$\;
	Set $\hat V_N$ $=$ $g$\;
	\For{$n$ $=$ $N-1,\ldots,0$}{
		Compute for $(z,a) \in \Gamma_{n} \times A$, 
		\beq
		\hat Q_n(z,a) & = &  f(z,a)+\sum_{\ell=1}^K p_\ell \widehat{V}_{n+1}\Big( \text{Proj}_{\Gamma_{n+1}} \big(F(z,a,e_\ell )\big) \Big),  \label{eq:estimateQvalueQknn}
		\enq
		where Proj$_{\Gamma_{n+1}}$ is the Euclidean projection over $\Gamma_{n+1}$\;
		\Comment{$\hat Q_n$ is the approximated $Q$-value\footnotemark{}  at time $n$}\\
		Compute the optimal control at time $n$ 
		\begin{equation}
		\label{eq:optimizeQknn}
		\hat A_n(z) \in \argmin_{a \in A}\big[ \hat Q_n(z,a)    \big], \quad \forall z \in \Gamma_{n};
		\end{equation}
		\Comment{use classical optimization algorithms of deterministic functions for this step}
		Set
		$
		\widehat{V}_{n}(z) \; = \;  \hat Q_n \big(z,\hat A_n(z) \big), \;\;\; \forall z \in \Gamma_n
		$\;
		\Comment{$\widehat{V}_{n}$ is the estimate by quantization of the value function}\\
	}
\end{algorithm}

\footnotetext{The $Q$-value at time $n$, denoted by $Q_n$, is defined as the function that takes the couple state-action $(x,a)$ as argument, and returns the expected optimal reward earned from time $n$ to time $N$ when the process $X$ is at state $x$ and action $a$ is chosen at time $n$; i.e. $Q_n: \R^d \times \R^q \backin  (x,a) \mapsto f(x,a) + \E_{n,x}^a[V_{n+1}(X_{n+1})]$.}

\begin{Remark}
{\rm	The estimate of the Q-value at time $n$ given by \eqref{eq:estimateQvalueQknn} is not continuous w.r.t. the control variable $a$, which might cause some stability issues when running Qknn, especially during the optimization procedure \eqref{eq:optimizeQknn}. We refer to  Section 3.2.2. in \cite{BalHurPha18} for a detailed presentation of an extension of Algorithm \ref{algo:Qknn} where the estimates of the $Q$ value function $Q_n$ is continuous w.r.t. the control variable.
}
\ep
\end{Remark}

\section{Numerical applications} \label{secnum}

\setcounter{equation}{0} \setcounter{Assumption}{0}
\setcounter{Theorem}{0} \setcounter{Proposition}{0}
\setcounter{Corollary}{0} \setcounter{Lemma}{0}
\setcounter{Definition}{0} \setcounter{Remark}{0}

In this section, we test the Neural-Networks-based algorithms presented in Section \ref{secalgo} on different examples. In high-dimension, we first took the same example as  already considered in \cite{Ehanjen17} so that we can directly compare our results to theirs, and take another example from linear quadratic control problem with explicit analytic solution that is served as reference value.
In low-dimension, we compared the results of our algorithms to the ones provided by Qknn, which has been introduced in Section \ref{secalgo} as an excellent benchmark for low-dimensional  control problems.

\subsection{A semilinear PDE}

\begin{figure}[H]
	\centering
	\includegraphics[width=0.9\linewidth]{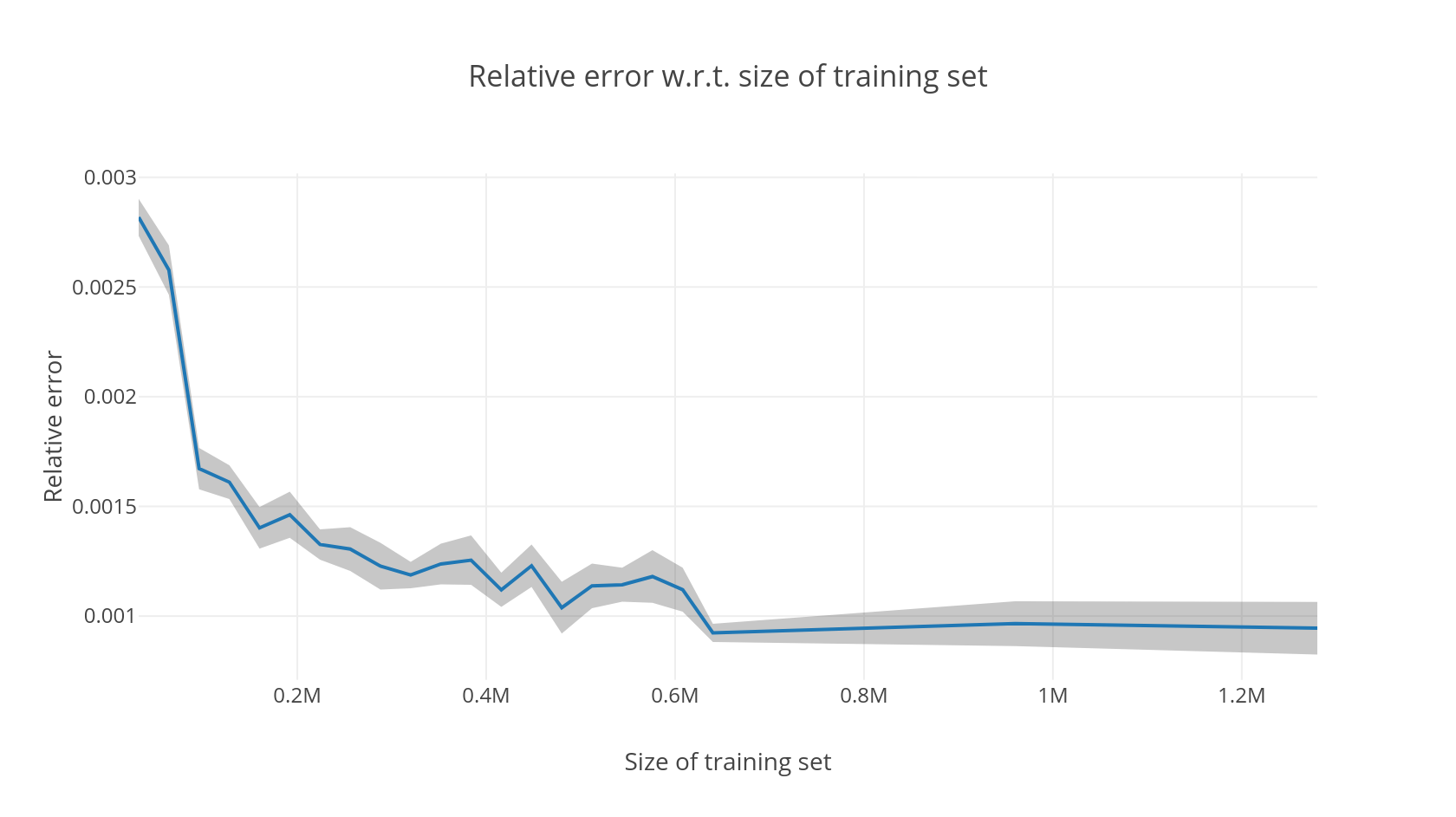}
	\vspace{-2mm}
	\caption{\small{Relative error of the Hybrid-Now estimate of the value function at time 0 w.r.t the number of mini-batches used to build the Hybrid-Now estimators of the optimal strategy. The value functions have been computed running three times a forward Monte Carlo with a sample of size 10,000, following the optimal strategy estimated by the Hybrid-Now algorithm.}}
	\label{fig:vfwrtmbisv2}
\end{figure}

We consider the following semilinear PDE with quadratic growth in the gradient:
\begin{equation} \label{HJBsemi}
\left\{
\begin{array}{rcl}
\Dt{v} +  \Delta_x v   \; -  \;   | D_x v |^2  &=& 0, \:\;\;\; (t,x) \in [0,T)\times\R^d, \\
v(T,x) & =& g(x), \;\;\; x \in \R^d. 
\end{array}
\right.
\end{equation}
By observing that for any $p$ $\in$ $\R^d$, -$|p|^2$ $=$ $\inf_{a \in \R^d} [ |a|^2 +  2 a.p]$, the PDE \reff{HJBsemi} can be written as a Hamilton-Jacobi-Bellman equation
\begin{equation} \label{HJBsemi2}
\left\{
\begin{array}{rcl}
\Dt{v} +  \Delta_x v \; +  \;   \inf_{a \in \R^d} \big[ |a|^2  + 2 a. D_x v ] &=& 0, \:\;\;\; (t,x) \in [0,T)\times\R^d, \\
v(T,x) & =& g(x), \;\;\; x \in \R^d,
\end{array}
\right.
\end{equation}
hence associated with the stochastic control problem
\beq \label{exconsto}
v(t,x) &=&  \inf_{ \alpha \in \Ac} \E \left[ \int_t^T |\alpha_s|^2 ds + g(X_T^{t,x,\alpha}) \right],  
\enq
where $X$ $=$ $X^{t,x,\alpha}$ is the controlled process governed by
\beqs
dX_s &=& 2 \alpha_s ds \; + \;  \sqrt{2} dW_s, \;\;\; t \leq s \leq T, \; X_t = x, 
\enqs
$W$ is a $d$-dimensional Brownian motion, and the control  process  $\alpha$ is valued in $A$ $=$ $\R^d$. The time discretization (with time step $h$ $=$ $T/N$) of the control problem \reff{exconsto} leads to the discrete-time control problem \reff{dynX}-\reff{costJ}-\reff{defcontrol} with
\beqs
X_{n+1}^\alpha &=& X_n^\alpha + 2 \alpha_n h + \sqrt{2 h} \eps_{n+1} \; =: \; F(X_n^\alpha,\alpha_n,\eps_{n+1}), \;\;\; n=0,\ldots,N-1,
\enqs
where $(\eps_n)_n$  is a sequence of i.i.d. random variables with law $\Nc(0,\I_d)$, and the cost functional
\beqs
J(\alpha) &=& \E \left[ \sum_{n=0}^{N-1} h |\alpha_n|^2 \; + \; g(X_N^\alpha) \right]. 
\enqs
On the other hand, it is known that an explicit solution to \reff{HJBsemi} (or equivalently \reff{HJBsemi2}) can be obtained via a Hopf-Cole transformation (see e.g. \cite{charic16}), and is given by
\beq
\label{eq:semiClosedForm}
v(t,x) &=&  -  \ln \Big( \E \Big[\exp\big(- g(x +  \sqrt{2} W_{T-t}) \big)\Big]  \Big), \;\;\; (t,x) \in [0,T]\times\R^d.  
\enq
We choose to run tests on two different examples that have already been considered in the literature:

\paragraph{Test 1} Some recent numerical results have been obtained in \cite{Ehanjen17} (see Section 4.3 in \cite{Ehanjen17}) when $T=1$ and $g(x)$ $=$ $\ln(\frac{1}{2}(1+|x|^2))$ in dimension $d$ $=$ $100$ (see Table 2 and Figure 3 in \cite{Ehanjen17}). Their method is based on neural network regression to solve the BSDE representation associated with the PDE \eqref{HJBsemi}, and provide estimates of  the value function at time 0 and state 0 for different values of a coefficient $\gamma$. We plotted the results of the Hybrid-Now algorithm in Figure \ref{fig:vfwrtmbisv2}.  Hybrid-Now took one hour to achieves a relative error of 0.11\%, using a 4-cores 3GHz intel Core i7 CPU. We want to highlight the fact that the algorithm presented in \cite{Ehanjen17} only needed 330 seconds to provide a relative error of 0.17\%. However, in our experience, it is difficult to reduce the relative error from 0.17\% to 0.11\% using their algorithm. Also, we believe that the computation time of our algorithm can easily be reduced; some ideas in this direction are discussed in Section \ref{secconclu}. The main trick that can be used is the transfer learning (also referred to as pre-training in the literature): we rely on the continuity of the value function and the optimal control w.r.t. time to claim that the value function and the optimal control at time $n$ are very close to the ones at time $n+1$. Hence, one can initialize the weights of the value function and optimal control at time $n$ with the optimal ones estimated at step $n+1$,  reduce the learning rate of the optimizer algorithm, and reduce the number of steps for the gradient descent algorithm. All this procedure really speeds up the learning of the value function and the optimal control, and insures stability of the estimates. Doing so, we were able to reduce the computation time from one hour to twenty minutes. 

\begin{figure}[H]
	\centering
	\includegraphics[width=.7\linewidth]{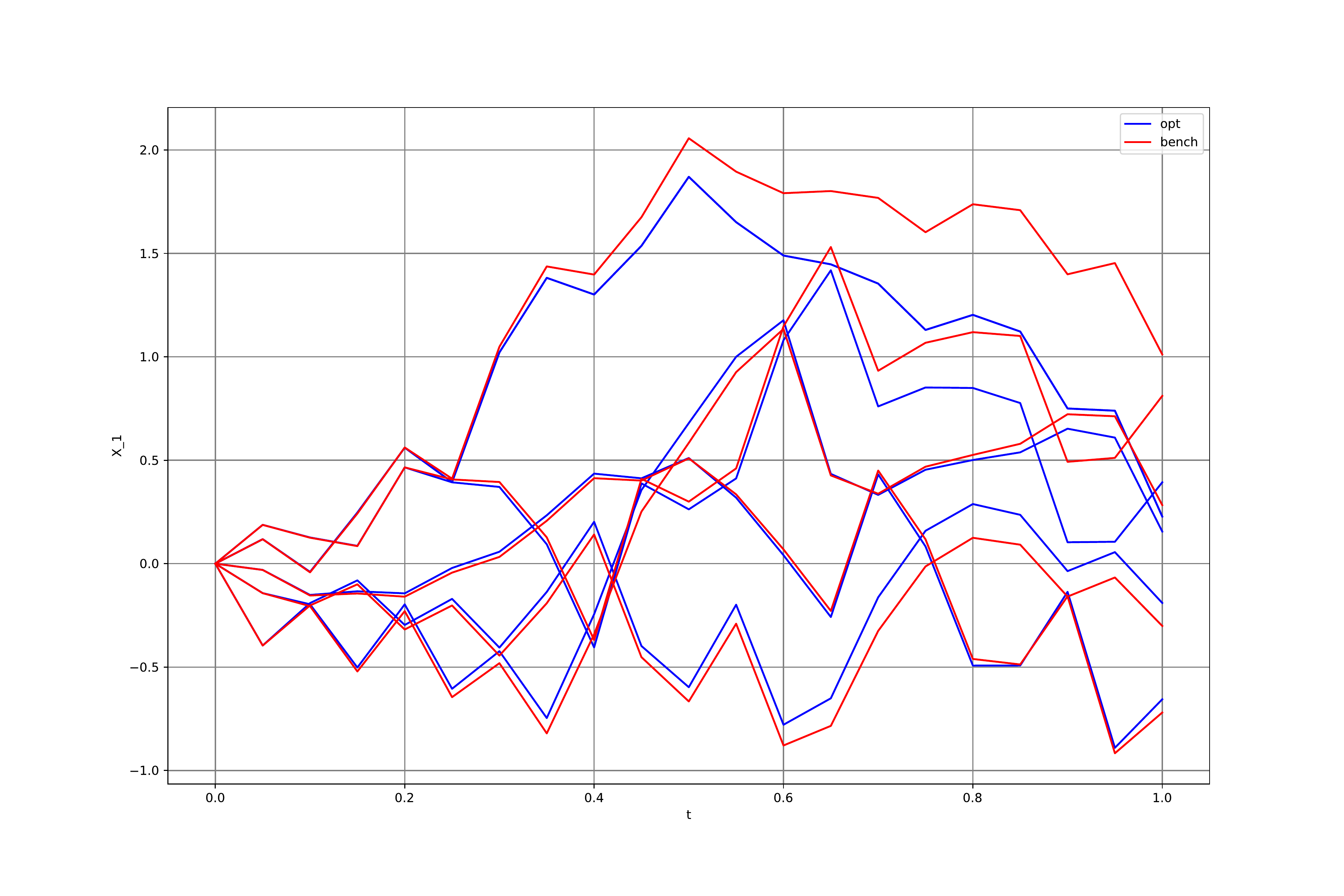}
	\vspace{-5mm}
	\caption{\small{Five forward simulations of the first component of $X$ w.r.t. time, when the agent follows the optimal strategy estimated by the Hybrid-Now (\bl{opt} in blue) and the naive strategy $\alpha=0$ (\red{bench} in red). We consider the problem in dimension $d$=2. Observe that the optimal strategy (estimated by Hybrid-Now) is to do nothing when the terminal time is far in order to avoid any running cost, i.e. $\alpha^{opt}=0$; and push $X$ toward 0 when the terminal time is close, in order to minimize the terminal cost.}}
	\label{fig:plotx1wrttpst1}
\end{figure}

We also considered the same problem in dimension $d$ $=$ $2$, for which we plotted the first component of $X$ w.r.t. time in Figure \ref{fig:plotx1wrttpst1}, for five different paths of the Brownian motion, where for each $\omega$, the agent follows either the naive ($\alpha$ $=$ $0$) or the Hybrid-Now strategy. One can see that both strategies are very similar when the terminal time is far; but the Hybrid-Now strategy clearly forces $X$ to get closer to 0 when the terminal time gets closer, in order to reduce the terminal cost.

\vspace{3mm}

Let us provide further implementation details on the algorithms presented in Test 1: 
\begin{itemize}
	\item As one can guess from the representation of $v$ in \eqref{exconsto}, it is probably optimal to drive the process $X$ around 0. Hence we decided to take $\mu_n:=(\frac{nT}{N})^{1/2}\mathcal{N}_d(0,I_d)$ as a training measure at time $n$ to learn the optimal strategy and value function at time $n$, for $n=0,\ldots,N-1$.
	\item We tested the algorithm with 1, 2 and 3 layers for the representation of the value function and the optimal control by neural networks, and noticed that the quality of the estimate significantly improves when using more than one layer, but does not vary significantly when considering more than 3 layers.
\end{itemize}

\paragraph{Test 2} Tests of the algorithms are proposed in dimension 1 with the terminal cost $g(x)= -x^\gamma \mathds{1}_{0 \leq x \leq 1} -\mathds{1}_{1 \leq x}$ and $\gamma \in (0,1)$. This problem was already considered in \cite{Richou10}, 
where the author proposed an algorithm based on a smart temporal discretization of the BSDE representation of the PDE \eqref{HJBsemi} in order to deal with the quadratic growth of the driver of the BSDE, and usual projection on basis functions techniques for the approximation of conditional expectations that appear in the dynamic programming equation associated with the BSDE. We refer to equations (13),(14),(15) in  \cite{ric11} for details on the proposed algorithm, and its Theorem 4.14 for the convergence result. Their estimates of the value function at time 0 and state 0, when $\gamma=1, 0.5, 0.1, 0$, are available in \cite{Richou10}, and have been reported in the column $Y\&R$ of Table \ref{t:semiLinearPDE}. Also, the exact values for the value function have been computed for these values of $\gamma$ by Monte Carlo using the closed-form formula \eqref{eq:semiClosedForm}, and are reported in the column \textit{Bench} of Table \ref{t:semiLinearPDE}.
	 Tests of the Hybrid-Now and Hybrid-LaterQ algorithms have been run, and the estimates of the value function at time 0 and state $x$ $=$ $0$ are reported in the Hybrid-Now and Hybrid-LaterQ columns. We also tested Qknn and reported its results in column Qknn. 
	 Note that Qknn is particularly well-suited to 1-dimensional control problems. In particular, it is not time-consuming since the dimension of the state space is $d$=1. Actually, it provides the fastest results, which is not surprising since the other algorithms need time to learn the optimal strategy and value function through gradient-descent method at each time step $n=0,\ldots,N-1$. Moreover, Table \ref{t:semiLinearPDE} reveals that Qknn is the most accurate algorithm on this example, probably because it uses local methods in space to estimate the conditional expectation that appears in the expression of the $Q$-value.

 \begin{table}[H]
	\centering
	\caption{\small{Value function at time 0 and state 0 w.r.t. $\gamma$, computed with the Y\&R, Hybrid-Now, Hybrid-Later and Qknn algorithms. Bench reports the MC estimates of the closed-form formula \eqref{eq:semiClosedForm}.}}
	\label{t:semiLinearPDE}
	\begin{tabular}{l|lllll}
		$\gamma$ & Y\&R      & Hybrid-LaterQ & Hybrid-Now   & Qknn  & Bench    \\ \hline
		1.0   & -0.402 & -0.456 & -0.460 & -0.461 & -0.464 \\
		0.5   & -0.466 & -0.495 & -0.507 & -0.508 & -0.509 \\
		0.1   & -0.573 & -0.572 & -0.579 & -0.581 & -0.586 \\
		0.0   & -0.620 & -1.000 & -1.000 & -1.000 & -1.000
	\end{tabular}
\end{table}

We end this paragraph by giving some implementation details for the different algorithms as part of  Test 2:
\begin{itemize}
\item\emph{Y\&R:} The algorithm Y\&R converged only when using a Lipschitz version of $g$. The following approximation was used to obtain the results in Table \ref{t:semiLinearPDE}:
 \[
 g_N(x)=
 \begin{cases}
 g(x) & \text{ if } x \not\in [0, N^{\frac{-1}{1-\gamma}}] \\
 -Nx & \text{ otherwise}.
 \end{cases}
 \]
 \item \emph{Hybrid-Now:}
We used $N=40$ time steps for the time-discretization of $[0,T]$. The value functions and optimal controls at time $n$ $=$ $0, \ldots,N-1$ are estimated using neural networks with 3 hidden layers and 10+5+5 neurons.
\item  \textit{Hybrid-LaterQ:}
We used $N=40$ time steps for the time-discretization of $[0,T]$. The value functions and optimal controls at time $n$ $=$ $0,\ldots,N-1$ are estimated using neural networks with 3 hidden layers containing 10+5+5 neurons; and 51 points for the quantization of the exogenous noise. 
\item  \textit{Qknn:}
 We used $N=40$ time steps for the time-discretization of $[0,T]$. We take 51 points to quantize the exogenous noise, $\eps_n \sim \mathcal{N}(0,1)$, for $n$ $=$ $0,\ldots, N$; and decided to use the 200 points of the optimal grid of $\mathcal{N}_2(0,1)$ for the state space discretization.
 \end{itemize}

The main conclusion regarding the results in this semilinear PDE problem is that Hybrid-Now provides better estimates of the solution to the PDE in dimension $d$=100 than the previous results available in \cite{Ehanjen17} but requires more time to do so. 
	
Hybrid-Now and Hybrid-Later provide better results than those available in \cite{ric11} to solve the PDE in dimension 2; but are outperformed by Qknn, which is arguably very accurate.

\subsection{A linear quadratic stochastic test case}

 We consider a  linear controlled process with dynamics in $\R^d$ according to 
 \beq \label{linXqua}
 dX_t &=& (BX_t + C \alpha_t) dt + \sum_{j=1}^p D_j \alpha_t dW^j_t, 
 \enq
 where $W^j$, $j$ $=$ $1,\ldots,p$, are independent real-Brownian motion,  the control process $\alpha$ $\in$ $\Ac$ is valued in $\R^m$, and the constant coefficients 
 $B$ $\in$ $\R^{d\times d}$, $C,D_j$ $\in$ $\R^{d\times m}$, $j$ $=$ $1,\ldots,p$.  The value function of the linear quadratic stochastic control problem is 
 \beqs
 v(t,x) &=& \inf_{\alpha\in\Ac} \E \Big[ \int_t^T (X_s^{t,x,\alpha}.Q X_s^{t,x,\alpha}  +  \lambda |\alpha_t|^2) dt +  X_T^{t,x,\alpha}. P X_T^{t,x,\alpha}  \Big], \;\; (t,x) \in [0,T]\times\R^d,
\enqs
 where $X^{t,x,\alpha}$ is the solution to \reff{linXqua} starting from $x$ at time $t$, given a control process $\alpha$ $\in$ $\Ac$,  $P,Q$ are nonnegative symmetric 
 $d\times d$ matrices, and $\lambda$ $>$ $0$.  
 The Bellman equation associated with this stochastic control problem is a fully nonlinear equation in the form
 \beqs
 \Dt{v} + x.Qx + \inf_{a\in\R} \big[ (Bx + C a).D_x v + a\trans\big( \lambda I_m + \sum_{j=1}^p \frac{D_j\trans D_x^2 v D_j}{2}\big) a \big]  
 &=& 0,  \; \mbox{ on } [0,T)\times\R^d, \\
 v(T,x) &=& x.Px, \;\;\; x \in \R^d, 
 \enqs
 and it is well-known, see e.g. \cite{yonzhou99}, that an explicit solution is given by 
 \begin{equation}
 v(t,x) = x. K(t) x,  \label{lq_vf}
 \end{equation}
where $K(t)$  is a nonnegative symmetric $d\times d$ matrix, solution to the Riccati equation
\beq \label{riccati}
\dot K + B\trans K + K B + Q - K C(\lambda I_m + \sum_{j=1}^p D_j\trans K D_j)^{-1} C\trans K 
&=& 0, \;\;\; K(T) = P, 
\enq
while an optimal feedback control is equal to
\begin{equation}
a^*(t,x)= - \big(\lambda I_m+ \sum_{j=1}^p D_j\trans K(t) D_j\big)^{-1} C\trans K(t) x,  
\;\;\;  (t,x) \in [0,T)\times\R^d. \label{lq_oc}
\end{equation}

We numerically solve this problem by considering a  time discretization (with time step $h$ $=$ $T/N$), which leads to the discrete-time control problem with dynamics 
\beqs
X_{n+1}^\alpha &=& X_n^\alpha + (B X_n^\alpha + C \alpha_n) h +  D \alpha_n \sqrt{h}  \eps_{n+1} \; =: \; F(X_n^\alpha,\alpha_n,\eps_{n+1}), \;\;\; n=0,\ldots,N-1,
\enqs
where $(\eps_n)_n$  is a sequence of i.i.d. random variables with law $\Nc(0,1)$, and  cost functional
\beqs
J(\alpha) &=& \E \left[ \sum_{n=0}^{N-1}  \big(X_n^\alpha.Q X_n^\alpha + \lambda |\alpha_n|^2)h  \; + \; X_N^\alpha. PX_N^\alpha \right]. 
\enqs
 
 For the numerical tests, we take $m$ $=$ $1$, $p$ $=$ $d$, and the following parameters: 
 \beqs
  T \; = \; 1, \; N=20,  & & B = I_d,  \;\;   C =  \mathds{1}_d, \; D_j = (0,\ldots,\underbrace{1}_{\text{j-th term}},\ldots,0)^\top, \; j=1,\ldots,p, \\
 & & Q \; = \; P \; =  I_d, \quad \lambda = 1,
\enqs
where we denote $\mathds{1}_d:=(\; \underbrace{1,\ldots,1}_{d\text{ times}} \;)^\top$.

\begin{figure}[H]
	\centering
	\includegraphics[width=1\linewidth]{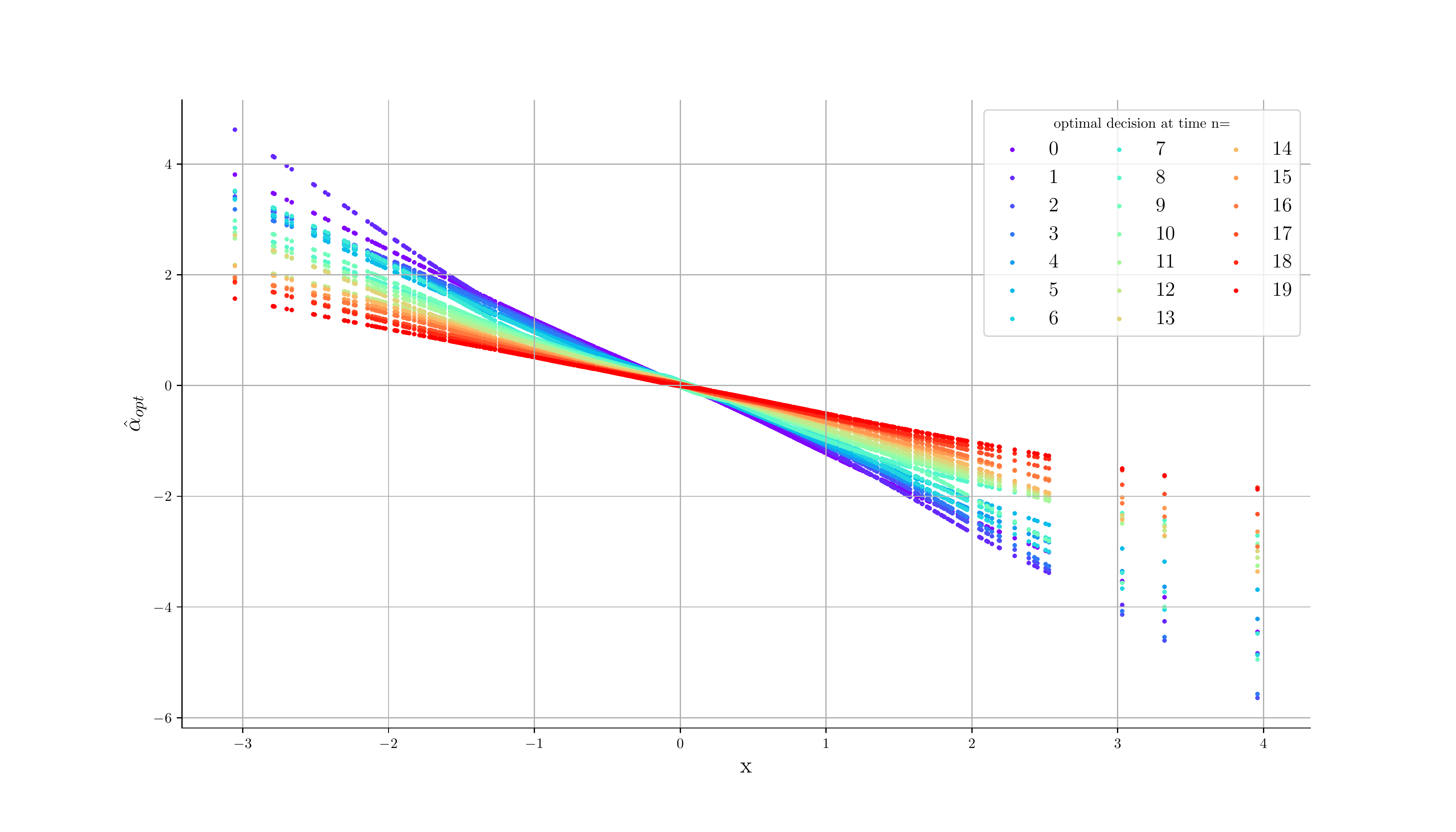}
	\caption{\small{Optimal decision estimated by Hybrid-Now at time $n=0,\ldots,N-1$. We took $d=1$, $N=20$. We observe that the estimates are linear, as expected given the closed-form formula \eqref{lq_oc} for the optimal control. }}
	\label{optimal_decisions_d1}
\end{figure}

\begin{figure}[H]
	\centering
	\includegraphics[width=1\linewidth]{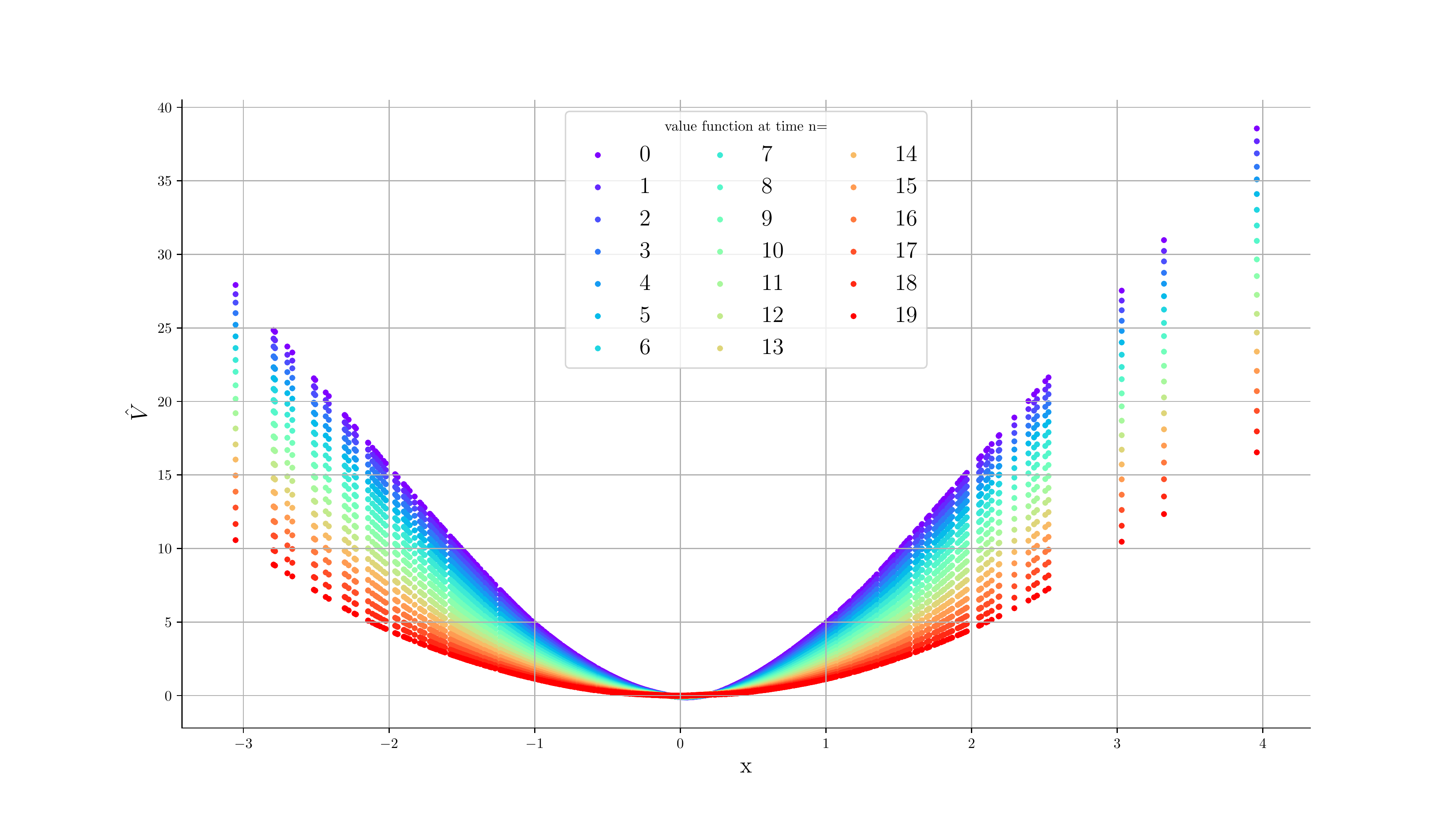}
	\caption{\small{Value function w.r.t. $x$, estimated by Hybrid-Now at time $n=0,\ldots,N-1$. We took $d=1$, $N=20$. We observe that the estimates are quadratic, as expected given the closed-form formula \eqref{lq_vf} for the value function.}}
	\label{vf_d1}
\end{figure}

\begin{figure}[H]
 	\centering
 	\includegraphics[width=1\linewidth]{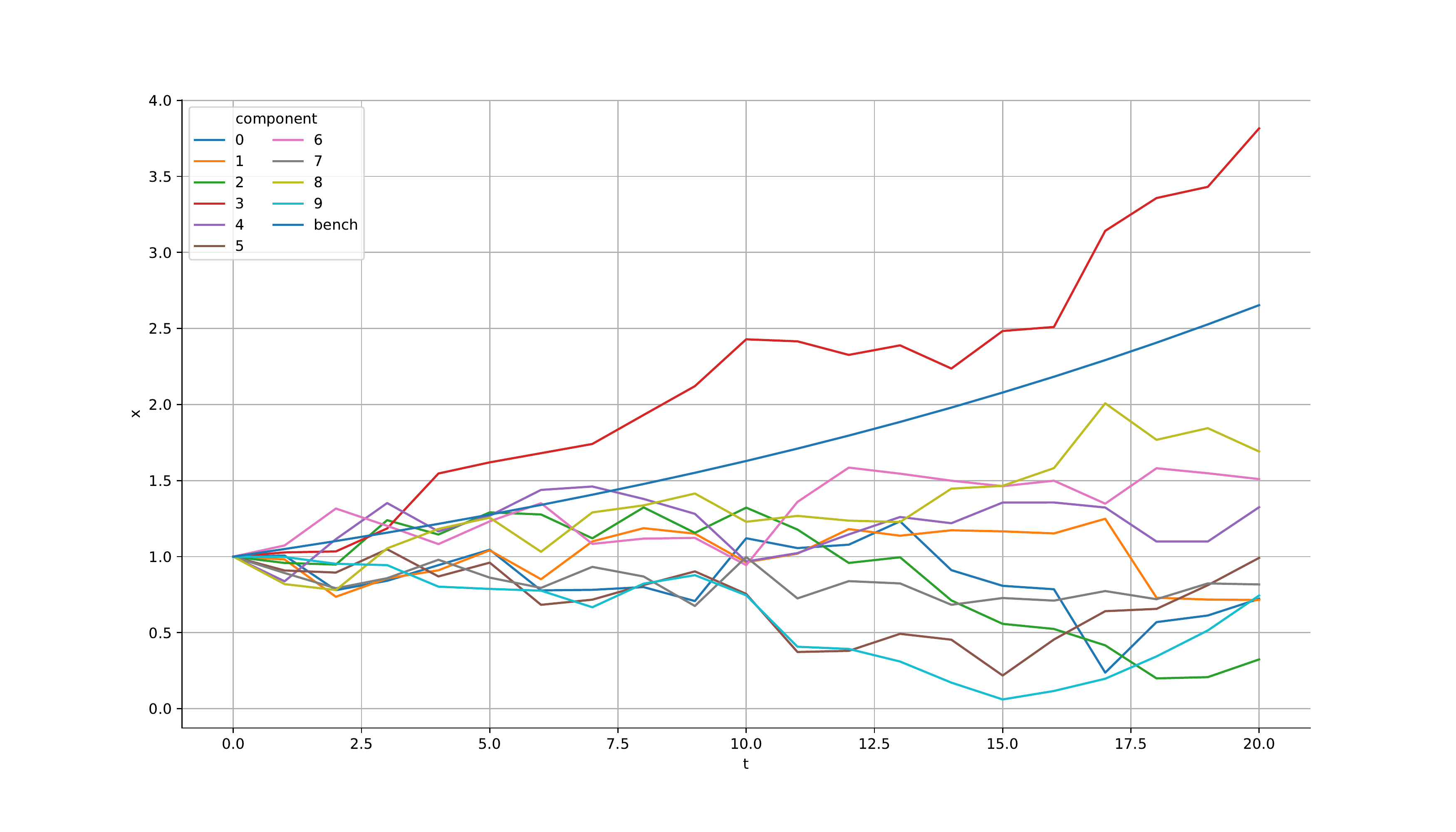}
 	\caption{\small{Forward simulation of $X$ w.r.t. time, when $X_0=1_d$ and $d=10$, driven optimally using Hybrid estimates. The first ten curves represent the ten components of $X$. The bench curve represents one of the identical component of $X$ when it is driven using the strategy $\alpha=0$. One can see that the optimal control tends to reduce the norm of each component of $X$.}}
 	\label{LQsimu_d10}
 \end{figure}

\paragraph{Numerical results} We implement our algorithms in dimension  $d=1,10,100$, and compare our solutions with the analytic solution via the Riccati equation \reff{riccati} solved by Matlab\footnote{We solved \eqref{riccati} with the Matlab method ode45.}. 
\begin{itemize}
	\item For $d$ $=$ $1$, we plotted the estimates of the optimal control at time $n=0,\ldots,N-1$ in Figure \ref{optimal_decisions_d1} and the value function in Figure \ref{vf_d1}. Observe that, as expected, the estimated optimal control is linear and the estimated value function is quadratic at each time. 
	\item For $d$ $=$ $10$, we reported in Table \ref{t:LQ_d10} the estimates of $v(0,X_0)$, computed by running forward simulations of $X$ using the estimated optimal strategy. ``Riccati'' is $v(0,X_0)$ computed by solving \eqref{riccati} with Matlab. We set the initial position to $X_0=\mathds{1}_d$. We also plotted in Figure \ref{LQsimu_d10} a forward simulation of the components of $X$ optimally controlled. Observe that NNContPI is more accurate than Hybrid-Now. Notice that the estimates provided by the algorithms are biased, which is due to the time discretization.
	\item For $d$ $=$ $100$, we reported in Table \ref{t:LQ_d100} the estimates of the value function, computed by running forward simulations of $X$ using the estimated optimal strategy. ``Riccati'' is $v(0,X_0)$ computed by solving \eqref{riccati} with Matlab. We set the initial position to $X_0=0.1 \mathds{1}_d$ and $X_0=0.5 \mathds{1}_d$. Once again, NNContPI is slightly more accurate than Hybrid-Now, and the estimates provided by the latter are biased due to the time discretization.
\end{itemize}

\emph{Implementation details:} We implemented Hybrid-Now and NNContPI using training sets from the distribution $\mu_n:=\mathcal{N}_d(0,1)$ for $n=0,\ldots,N-1$. We represented the value function and optimal control at time $n$, $n=0,\ldots,N-1$ using two hidden layers with d+20 and d+10 neurons, and 1 neuron for the output layers. We used Elu as activation function for the hidden layers, and identity for the output layer. \\


\newpage

\emph{Comments on the algorithms:}  Hybrid-Now behaved similarly  as for the SemiLinear PDE example, and we can make the same remarks. 
NNContPI is much slower than Hybrid-Now, because the data have to go through the $N-n-1$ neural networks that represent the optimal controls at time $n+1,\ldots,N-1$, 
in order to estimate the optimal control at time $n$. 

\begin{table}[H]
	\centering
	\caption{\small{Estimate of $v(0,X_0)$ obtained by forward simulation of the process controlled by the optimal strategy estimated by Hybrid-Now and NNContPI. ``Riccati'' is $v(0,X_0)$ computed by solving \eqref{riccati} with Matlab. We took $d=10$, and $X_0=\mathds{1}_d$. Mean and standard deviation are computed on 10 sets of 10,000 simulations each.}}
	\label{t:LQ_d10}
	\begin{tabular}{c|c|c}
		&	Mean   &   std    \\ \hline
		Hybrid-Now &	56.0 & 0.6\\
		NNContPI & 54.3 & 0.1 \\
		Riccati & 57.1 & -
	\end{tabular}
\end{table}

\begin{table}[H]
		\caption{\small{Estimate of $v(0,X_0)$ obtained by forward simulation of the process controlled by the optimal strategy estimated by Hybrid-Now and NNContPI. ``Riccati'' is $v(0,X_0)$ computed by solving \eqref{riccati} with Matlab. We took $d=100$, and initial position $X_0=0.5 \mathds{1}_d$ and $X_0=0.1 \mathds{1}_d$. Mean and standard deviation are computed on 10 sets of 10,000 simulations each.}}
		\label{t:LQ_d100}
	\begin{minipage}{.5\linewidth}
		\centering
	\begin{tabular}{c|c|c}
	&	Mean   &   std    \\ \hline
	Hybrid-Now &5.7 & 7e-3\\
	NNContPI & 5.4 & 7e-3 \\
	Riccati & 5.7 & - \\ \hline
\end{tabular}
\vspace{1mm}
\caption*{Case $X_0=0.1 \mathds{1}_d$}
	\end{minipage}
	\begin{minipage}{.5\linewidth}
	\centering
	\begin{tabular}{c|c|c}
	&	Mean   &   std    \\ \hline
	Hybrid-Now &	137.1 & 1.3e-1 \\
	NNContPI & 137.4& 1.4e-1 \\
	Riccati & 142.7 & - 	\\ \hline
\end{tabular}
\vspace{1mm}
\caption*{Case $X_0=0.5\mathds{1}_d$}
\end{minipage}
\end{table}

\subsection{Option hedging} 
  
Our third example comes from a classical hedging problem in finance. We consider an investor who trades in $q$ stocks with (positive) price process  $(P_n)_n$, and we denote by $(\alpha_n)$ valued in $\A$ $\subset$ $\R^q$ the amount held in these assets over the period $(n,n+1]$. We assume for simplicity that the price of the riskless asset is constant equal to $1$ (zero interest rate). It is convenient to introduce the return process as: $R_{n+1}$ $=$  ${\rm diag}(P_n)^{-1}(P_{n+1}-P_n)$, $n$ $=$ $0,\ldots,N-1$, so that the self-financed wealth process of the investor with a portfolio strategy $\alpha$, and starting from some capital $w_0$, is governed by 
\beqs
\Wc_{n+1}^\alpha &=& \Wc_n^\alpha + \alpha_n. R_{n+1}, \;\;\; n = 0,\ldots,N-1, \;\; \Wc_0^\alpha = w_0. 
\enqs
Given an option payoff $h(P_N)$, the objective of the agent is to minimize over her portfolio strategies $\alpha$ her expected square replication error
\beqs
V_0 &=& \inf_{\alpha\in\Ac} \E \Big[ \ell\big( h(P_N) - \Wc_N^\alpha \big) \Big], 
\enqs
where $\ell$ is a convex function on $\R$.  
Assuming that the returns $R_n$, $n$ $=$ $1,\ldots,N$ are i.i.d, we are in a $(q+1)$-dimensional framework of Section \ref{secintro} with 
$X^\alpha$ $=$ $(\Wc^\alpha,P)$ with $\eps_{n}$ $=$ $R_n$ valued in $E$ $\subset$ $\R^q$, with the dynamics function
\beqs
F(w,p,a,r) &=& \left\{ \begin{array}{c}
w + a.r \\
p + {\rm diag}(p)r,
\end{array}
\right. \;\;\;\;\;   x= (w,p) \in \R\times\R^q, \; a \in \R^q, \; r \in E, 
\enqs
the running cost function $f$ $=$ $0$ and the terminal cost $g(w,p)$ $=$ $\ell(h(p) - w)$.   
We test our algorithm in the case of a square loss function, i.e. $\ell(w)$ $=$ $w^2$, and when there is no portfolio constraints $\A$ $=$ $\R^q$, and compare our numerical results with 
the explicit solution derived in \cite{beretal01}:  denote by  $\nu(dr)$ the distribution of $R_n$, by $\bar\nu$ $=$ $\E[R_n]$ $=$ $\int r \nu(dr)$ its mean, and by 
$\bar M_2$ $=$ $\E[R_n R_n\trans]$ assumed to be invertible;  we then   have 
\beqs
V_n(w,p) &=& K_n w^2 -  2 Z_n(p) w + C_n(p)
\enqs
where the functions $K_n$ $>$ $0$, $Z_n(p)$ and $C_n(p)$ are  given in backward induction, starting from the terminal condition  
\beqs
K_N \;  = \; 1, \;\; Z_N(p) \; = \;  h(p), \;\; C_N(p) \; = \; h^2(p),
\enqs
and for $n$ $=$ $N-1,\ldots,0$,  by 
\begin{eqnarray*}
K_n &=& K_{n+1}  \big( 1 -  \bar\nu\trans \bar M_2^{-1}\bar\nu \big),   \\
Z_n(p) &=& \int Z_{n+1}(p + {\rm diag}(p)r) \nu(dr) -   \bar\nu\trans \bar M_2^{-1} \int Z_{n+1}(p + {\rm diag}(p)r) r \nu(dr), \\
C_n(p) &=& \int C_{n+1}(p + {\rm diag}(p)r) \nu(dr)  \\
& & \;\;\;\;\;\;\; - \;   \frac{1}{K_{n+1}}  \Big( \int  Z_{n+1}(p + {\rm diag}(p)r) r \nu(dr) \Big)\trans \bar M_2^{-1}    \Big( \int  Z_{n+1}(p + {\rm diag}(p)r) r \nu(dr) \Big),  
\end{eqnarray*}
so that $V_0$ $=$ $K_0 w_0^2 - 2 Z_0(p_0)w_0 + C_0(p_0)$, where $p_0$ is the initial stock price.  Moreover, the optimal portfolio strategy is given in feedback form by 
$\alpha_n^*$ $=$ $a^*_n(\Wc_n^*,P_n)$,  where $a^*_n(w,s)$ is the function
\[
a_n^*(w,p) =   \bar M_2^{-1}  \left[ \frac{\int Z_{n+1}(p + {\rm diag}(p)r) r \nu(dr)}{K_{n+1}}   -  \bar\nu w \right],
\]
and $\Wc^*$ is the optimal wealth associated with $\alpha^*$, i.e., $\Wc_n^*$ $=$ $\Wc_n^{\alpha^*}$.  Moreover, the initial capital $w_0^*$ that minimizes $V_0$ $=$ $V_0(w_0,p_0)$, and called (quadratic) hedging price is given by
\beqs
w_0^* &=& \frac{Z_0(p_0)}{K_0}. 
\enqs

 \vspace{2mm}

\paragraph{Test} 
Take $N=6$, and consider one asset $q$ $=$ $1$ with returns  modeled by a trinomial tree: 
 \beqs
 \nu(dr) &=&  \pi_+ \delta_{r_+} + \pi_0 \delta_0 + \pi_- \delta_{r_-}, \;\;\; \pi_0 + \pi_+  +  \pi_-  = 1, 
 \enqs
 with $r_+$ $=$ $5\%$, $r_-$ $=$ $-5\%$, $\pi_+$ $=$ $60\%$, $\pi_-$ $=$ $30\%$. Take $p_0$ $=$ $100$, and 
consider the call option $h(p)$ $=$ $(p-\kappa)_+$ with $\kappa$ $=$ $100$. The price of this option  is defined as the initial value of the portfolio that minimizes the terminal quadratic loss of the agent when the latter follows the optimal strategy associated with the initial value of the portfolio. In this test, we want to determine the price of the call and the associated optimal strategy using different algorithms.

\begin{Remark}
	{\rm 	The option hedging problem is linear-quadratic, hence belongs to the class of problems where the agent has ansatzes on the optimal control and the value function. Indeed, we expect here the optimal control to be affine w.r.t. $w$ and the value function to be quadratic w.r.t. $w$. For these kind of problems, the algorithms presented  in Section \ref{secalgo} can easily be adapted so that the expressions of the estimators satisfy the ansatzes. See \eqref{class::ann} and \eqref{class::F} for the option hedging problem. 
	}
	\ep
\end{Remark}
\vspace{2mm}

\begin{figure}
	\centering
	\makebox[\textwidth][c]{\includegraphics[width=1\linewidth]{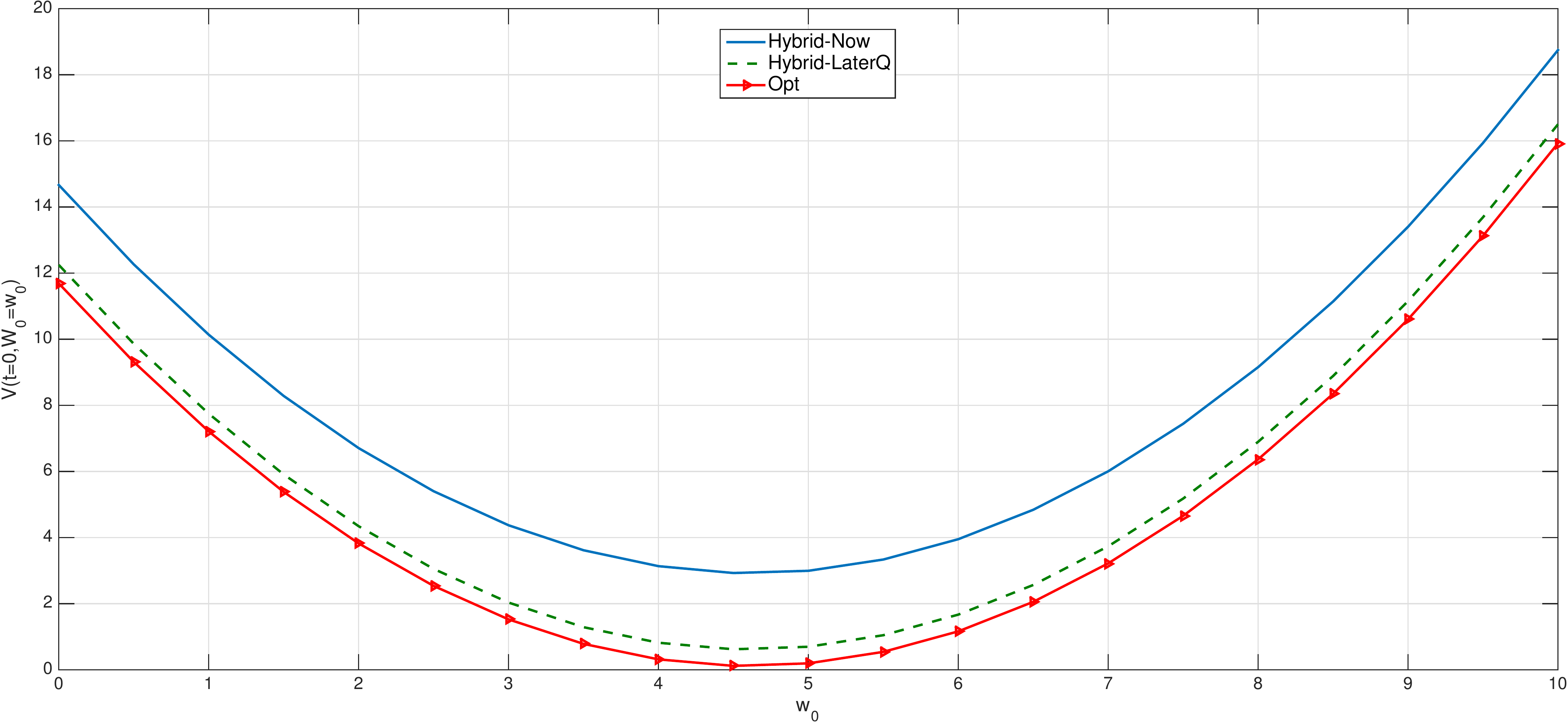}}
	\caption{Estimates of the value function at time 0 w.r.t. $w_0$ using Hybrid-Now (blue line) or Hybrid-LaterQ (green dashes). We draw the value function in red for comparison. One can observe that all the algorithms estimate the price to be 4.5, but Hybrid-LaterQ is better than Hybrid-Now at reducing the quadratic risk.}
	\label{fig:OHvf0}
\end{figure}

\paragraph{Numerical results}
In Figure \ref{fig:OHvf0}, we plot the value function at time 0 w.r.t $w_0$, the initial value of the portfolio, when the agent follows the theoretical optimal strategy (benchmark), and the optimal strategy estimated by the Hybrid-Now or Hybrid-LaterQ algorithms. We perform forward Monte Carlo using 10,000 samples to approximate the lower bound of the value function at time 0 (see \cite{PHL17} for details on how to get an approximation of the upper-bound of the value function via duality). One can observe that while all the algorithms give a call option price approximately equal to 4.5, Hybrid-LaterQ clearly provides a better strategy than Hybrid-Now to reduce the quadratic risk of the terminal loss. 

\begin{figure}
	\centering
	\includegraphics[width=.7\linewidth]{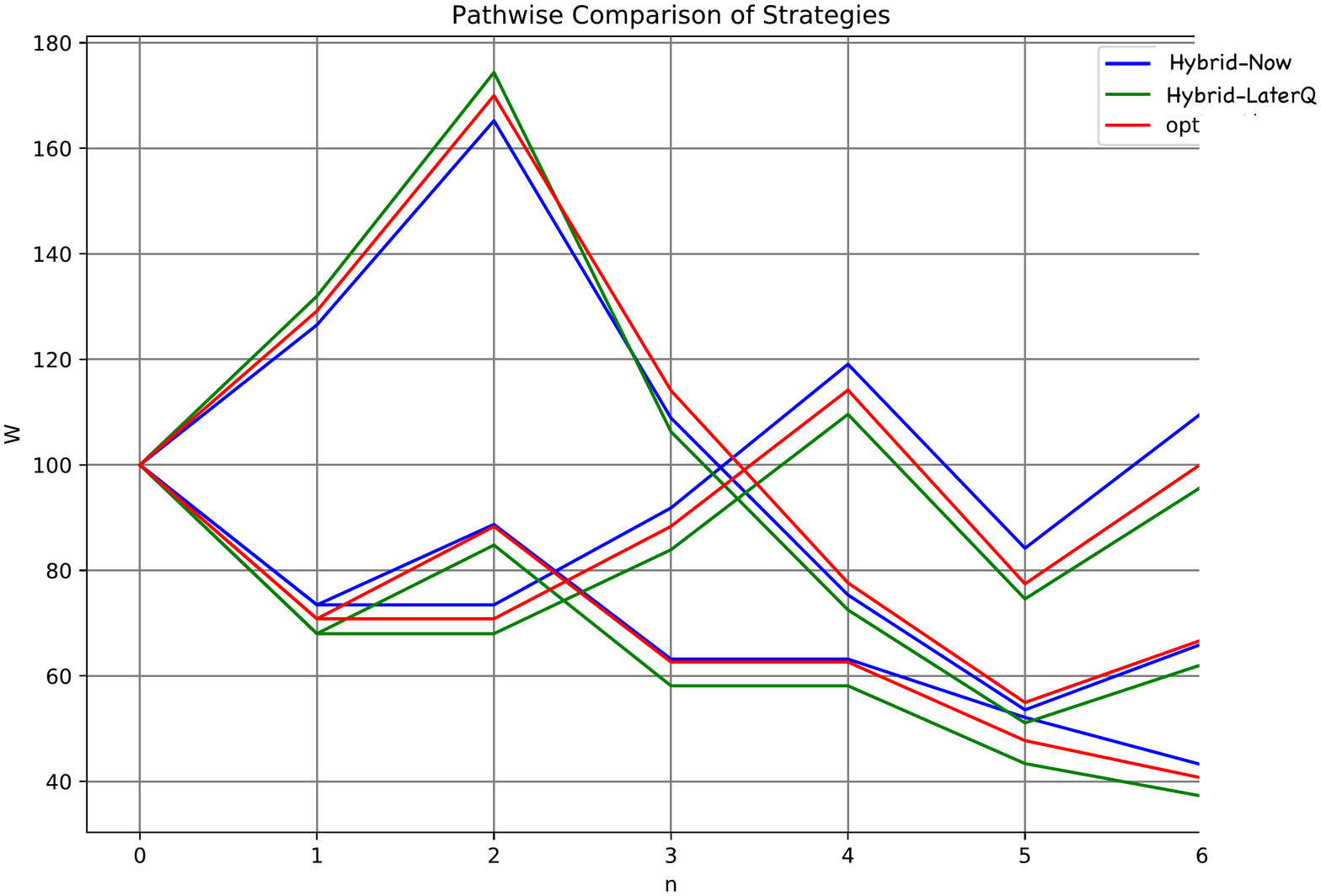}
	\caption{\small{Three simulations of the agent's wealth w.r.t. time $n$ when, for each $\omega$, the latter follows the theoretical optimal strategy (red), the estimated one using Hybrid-Now (blue) and the one using Hybrid-LaterQ (green). We took $w_0=100$. Observe that the process is driven similarly to the optimally controlled process, when the agent follows the estimated optimal strategy using Hybrid-LaterQ or Hybrid-Now.}}
	\label{fig:patwhisecomparison2}
\end{figure}

We plot in Figure \ref{fig:patwhisecomparison2} three different paths of the value of the portfolio w.r.t the time $n$,  when the agent follows either the theoretical optimal strategy (red), or the estimated one using Hybrid-Now (blue) or Hybrid-LaterQ (green). We set $w_0=100$ for these simulations. 

\vspace{2mm}

\paragraph{Comments on Hybrid-Now and Hybrid-LaterQ}
The Option Hedging problem belongs to the class of linear-quadratic control problems for which we expect the optimal control to be affine w.r.t. $w$ and the value function to be quadratic w.r.t. $w$. It is then natural to consider the following classes of controls $\mathcal{A}_M$ and functions $\mathcal{F}_M$ to properly approximate the optimal controls and the values functions at time $n$=$0,\ldots,N-1$:
\begin{equation}
\label{class::ann}
\mathcal{A}_M:= \left\{ (w,p) \mapsto A(x;\beta) \cdot  \big(1,w\big)^\intercal  ; \quad \beta \in \mathbb{R}^p \right\},
\end{equation}
\begin{equation}
\label{class::F}
\mathcal{F}_M:= \left\{ (w,p) \mapsto \Phi(x;\theta) \cdot  \big(1,w,w^2\big)^\intercal  ; \quad \theta \in \mathbb{R}^p \right\},
\end{equation}
where $\beta$ describes the parameters (weights+bias) associated with the neural network $A$ and $\theta$ describes those associated with the neural network $\Phi$. The notation $^\intercal$ stands for the transposition, and $\cdot$ for the inner product. Note that there are 2 (resp. 3) neurons in the output layer of $A$ (resp. $\Phi$), so that the inner product is well-defined in \eqref{class::F} and \eqref{class::ann}.

\subsection{Valuation of energy storage}\label{subsectionEnergyStorage}

We present a discrete-time version of the energy storage valuation problem studied in \cite{carlud07}. We consider a commodity (gas) that has to be stored in a cave, e.g. salt domes or aquifers. The manager of such a cave aims to maximize the real options value by optimizing over a finite horizon $N$ the dynamic decisions to inject or withdraw gas as time and market conditions evolve. We denote by 
$(P_n)$ the gas price, which is an exogenous real-valued Markov process 
modeled by the following mean-reverting process: 
\beq \label{dynP}
P_{n+1} &=& \bar p (1-\beta)  + \beta P_n + \xi_{n+1},
\enq
where $\beta$ $<$ $1$, and $\bar p$ $>$ $0$ is the stationary value of the gas price.  The current  inventory in the gas storage is denoted by $(C_n^\alpha)_n$ and depends on the manager's decisions represented by a control process $\alpha$ $=$ $(\alpha_n)$ valued in 
$\{-1,0,1\}$:  $\alpha_n$ $=$ $1$ (resp. $-1$) means that she injects (resp. withdraws) gas with an injection (resp. withdrawal) rate $a_{in}(C_n^\alpha)$ (resp. $a_{out}(C_n^\alpha)$) 
requiring (causing)  a purchase  (resp. sale) of $b_{in}(C_n^\alpha)$ $\geq$ $a_{in}(C_n^\alpha)$ (resp.  $b_{out}(C_n^\alpha)$ $\leq$ $a_{out}(C_n^\alpha)$), and $\alpha_n$ $=$ $0$ means that she is doing nothing. The difference between  $b_{in}$ and $a_{in}$ (resp. $b_{out}$ and $a_{out}$) indicates gas loss during injection/withdrawal.   
The evolution of the inventory is then governed by 
\beq \label{dynC}
C_{n+1}^\alpha &=& C_n^\alpha  +  h(C_n^\alpha,\alpha_n), \;\;\; n =0,\ldots,N-1, \; C_0^\alpha \; = \; c_0,
\enq 
where we set
\beqs
h(c,a) &=& \left\{ \begin{array}{cl}
a_{in}(c) & \mbox{ for } \; a =1 \\
0 &  \mbox{ for } \; a = 0 \\
- a_{out}(c) & \mbox{ for } \; a = - 1,
\end{array}
\right.
\enqs
 and we have the physical inventory constraint:
\beqs
C_n^\alpha  & \in & [C_{min},C_{max}], \;\;\; n =0,\ldots,N. 
\enqs 
The running gain of the manager at time $n$ is $f(P_n,C_n^\alpha,\alpha_t)$ given by 
\beqs
f(p,c,a) &=& \left\{ \begin{array}{cl}
- b_{in}(c) p  - K_{1}(c)    & \mbox{ for } \; a =1 \\
- K_0(c) &  \mbox{ for } \; a = 0 \\
b_{out}(c) p - K_{-1}(c) & \mbox{ for } \; a = -1,
\end{array}
\right.
\enqs
and  $K_{i}(c)$ represents the storage cost in each regime $i$ $=$ $-1,0,1$. The problem of the manager is then to maximize over $\alpha$ the expected total profit
\beq
\label{eq:rewardVES}
J(\alpha) &=& \E \left[ \sum_{n=0}^{N-1} f(P_n,C_n^\alpha,\alpha_n) + g(P_N,C_N^\alpha)\right],
\enq
where a common choice for the terminal condition is
\beqs
g(p,c) &=& -   \mu p(c_0 - c)_+, 
\enqs
which penalizes for having less gas than originally, and makes this penalty proportional to the current price of gas ($\mu$ $>$ $0$). 
We are then in the $2$-dimensional framework of Section \ref{secintro} with $X^\alpha$ $=$ $(P,C^\alpha)$, and the set of admissible controls in the dynamic programming loop is given by:
 \beqs
 A_n(c) = \big\{ a \in \{-1,0,1\}: c + h(c,a) \in [C_{min},C_{max}], \;c \in [C_{min},C_{max}] \big\}, \;\; n =0,\ldots,N-1. 
 \enqs

\vspace{3mm}

\paragraph{Test} 
We fixed the parameters as follows, to run our numerical tests: 
\beqs
a_{in}(c) \; = \; b_{in}(c) \; = \; 0.06, & & a_{out}(c) \; = \; b_{out}(c) \; = \; 0.25 \\
K_i(c) \; = \; 0.01 c
\enqs
$C_{max}$ $=$ $8$,  $C_{min}$ $=$ $0$, $c_0$ $=$ $4$,  $\bar p$ $=$ $5$, $\beta$ $=$ $0.5$, $\xi_{n+1}$ $\leadsto$ ${\cal N}(0,\sigma^2)$ with $\sigma^2$ $=$ $0.05$,  and 
$\mu$ $=$ $2$ in the terminal penalty function, $N$ $=$$30$.    

\vspace{2mm}

\begin{figure}[H]
	\centering
	\makebox[\textwidth][c]{\includegraphics[width=1.2\linewidth]{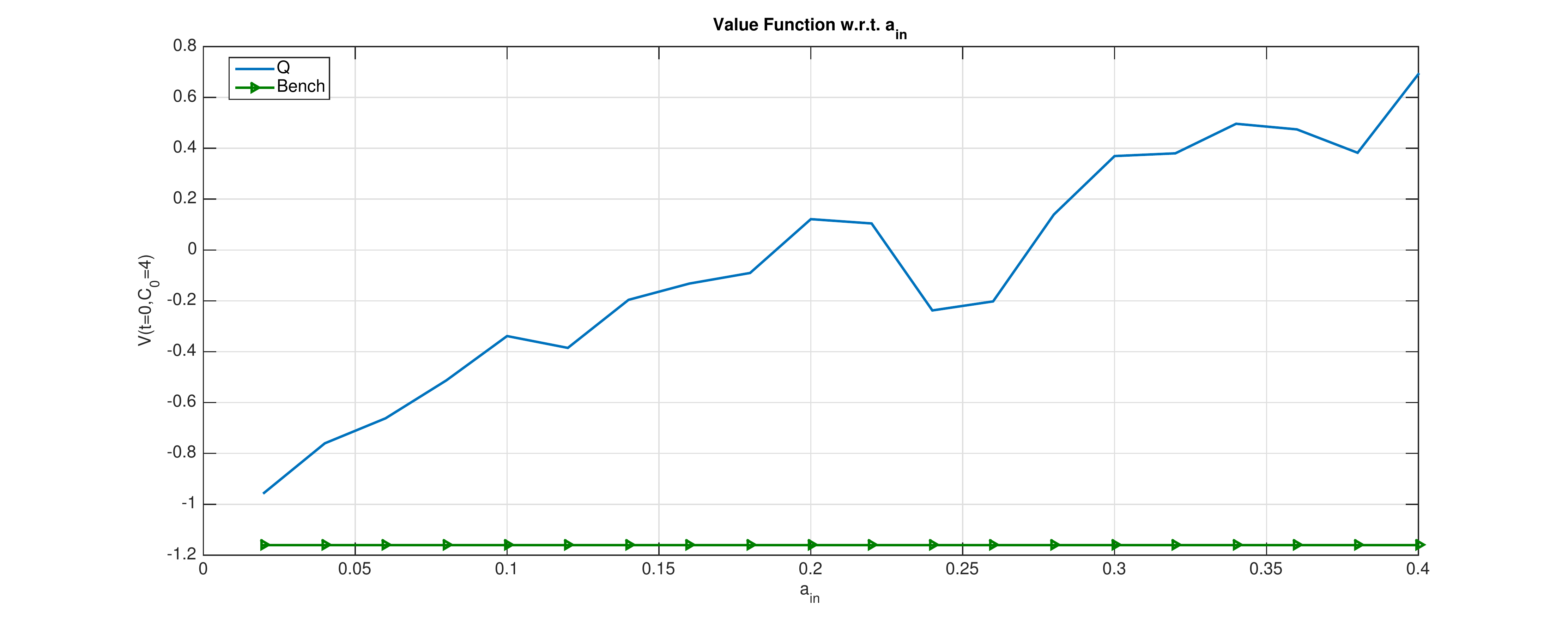}}
	\caption{\small{Estimate of the value function at time 0 w.r.t. $a_{in}$, when the agent follows the optimal strategy estimated by Qknn, by running a forward Monte Carlo with a sample of size 100,000 (blue). We also plotted the cost functional associated with the naive passive strategy $\alpha$ $=$ $0$ (Bench). See that for small values of $a_{in}$ such as 0.06, doing nothing is a reasonable strategy. Observe also that the value function is not monotone w.r.t. $a_{in}$ which is due to the dynamics of $C$ \eqref{dynC}.}}
	\label{fig:VESvf0}
\end{figure}

\paragraph{Numerical results} We plotted in Figure \ref{fig:VESvf0} the estimates of the value function at time 0 w.r.t. $a_{in}$ using Qknn, as well as the reward function \eqref{eq:rewardVES} associated with the naive do-nothing strategy $\alpha$ $=$ $0$ (see Bench in figure \ref{fig:VESvf0}). As expected, the naive strategy performs well when $a_{in}$ is small compared to $a_{out}$, since, in this case, it takes time to fill the cave, so that the agent is likely to do nothing in order to avoid any penalization at terminal time. When $a_{in}$ is of the same order as $a_{out}$, it is easy to fill up and empty the cave, so the agent has more freedom to buy and sell gas in the market without worrying about the terminal cost. 
Observe that the value function is not monotone,  due to the fact that the $C$ component in the state space takes its value in a bounded and discrete set (see \eqref{dynC}).

Table \ref{t:VES} provides the estimates of the value function using the \cmt{{\color{magenta}  AlgINQVI, AlgINQPI,}} ClassifPI, Hybrid-Now and Qknn algorithms.
Observe first that the estimates provided by Qknn are larger than those provided by the other algorithms, meaning that Qknn outperforms the other algorithms. The second best algorithm is ClassifPI, while Hybrid-Now performs poorly and clearly suffers from instability, due to the discontinuity of the running rewards w.r.t. the control variable. 

\begin{table}[H]
	\centering
	\caption{\small{$ V(0,P_0,C_0)$ estimates for different values of $a_{in}$, using the optimal strategy provided by the ClassifPI , Hybrid-Now  and Qknn algorithms, with $a_{out}$ $=$ $0.25$, $P_0$ $=$ $4$ and $C_0$ $=$ $4$.}}
	\label{t:VES}
	\begin{tabular}{l|llll}
	$a_{in}$  &Hybrid-Now& ClassifPI & Qknn  & $ \alpha=0$ \\ \hline
	0.06 & -0.99 & -0.71 & -0.66 & -1.20\\
	0.10 & -0.70 & -0.38 & -0.34 & -1.20\\
	0.20 & -0.21 &\ 0.01 &\ 0.12 & -1.20\\
	0.30 & -0.10 &\ 0.37 &\ 0.37 & -1.20\\
	0.40 &\ 0.10 &\ 0.51 &\ 0.69 & -1.20\\
\end{tabular}
\end{table}

\begin{figure}
	\begin{subfigure}{.5\linewidth}\centering
		{
			\includegraphics[width=1.2\linewidth]{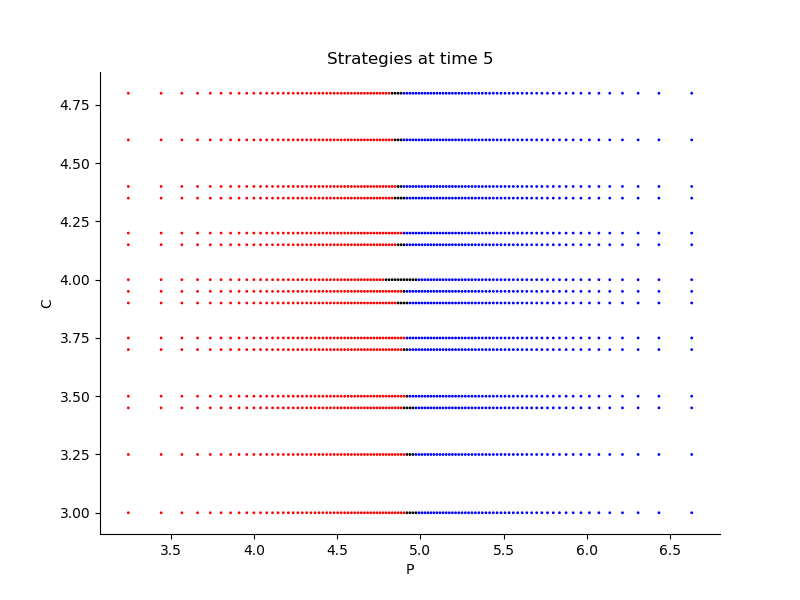}
		}
		\vspace*{-1cm}
		\caption*{}%
	\end{subfigure}%
	\begin{subfigure}{.5\linewidth}\centering
		{	
			\includegraphics[width=1.2\linewidth]{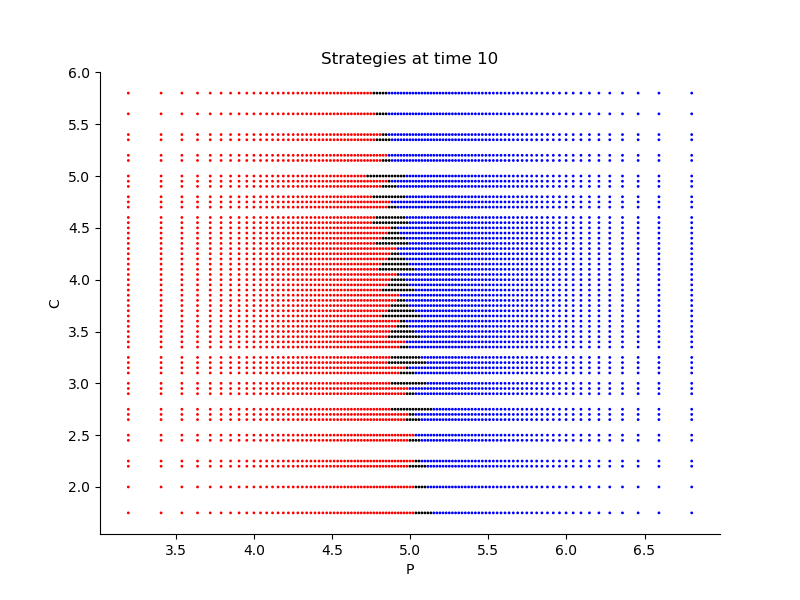}
		}
		\vspace*{-.8cm}
		\caption*{}%
	\end{subfigure}
	\begin{subfigure}{.5\linewidth}\centering
		{	
			\includegraphics[width=1.2\linewidth]{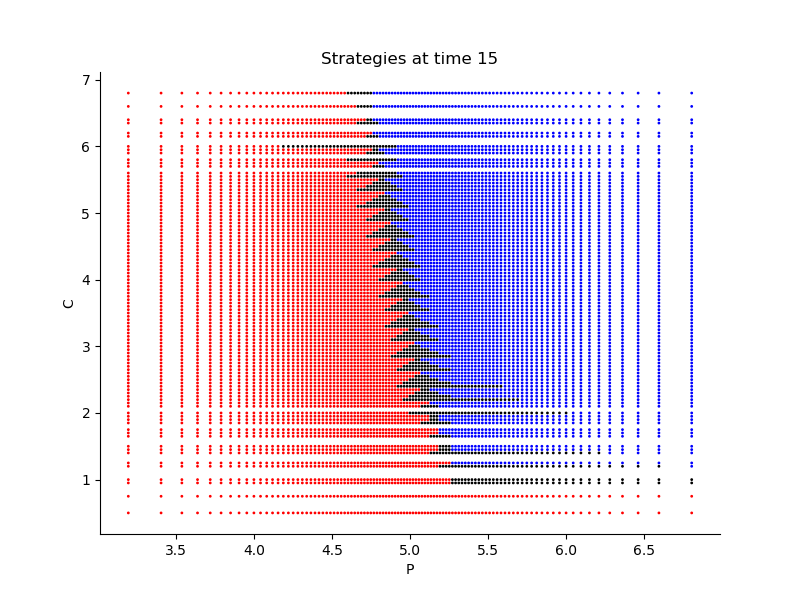}
		}
		\vspace*{-.8cm}
		\caption*{}%
	\end{subfigure}
	\begin{subfigure}{.5\linewidth}\centering
		{	
			\includegraphics[width=1.2\linewidth]{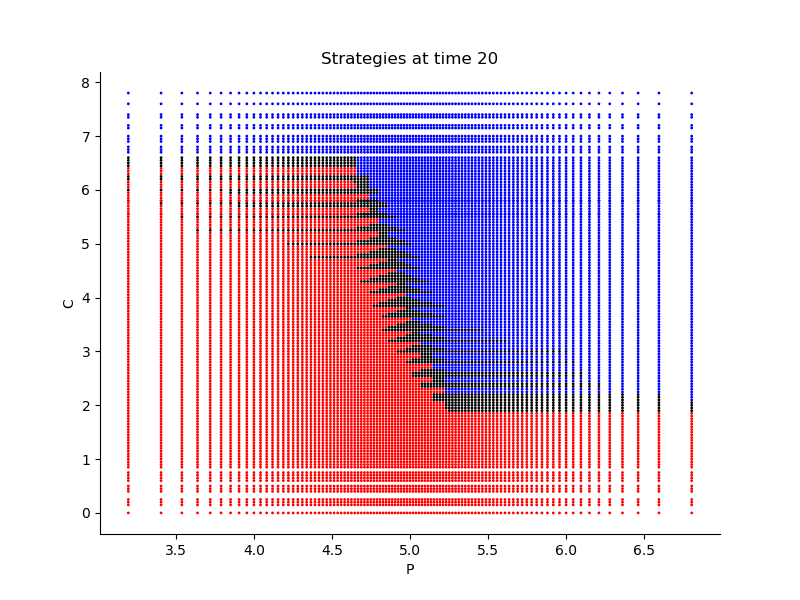}
		}
		\vspace*{-.8cm}
		\caption*{}%
	\end{subfigure}
	\begin{subfigure}{.5\linewidth}\centering
	{	
		\includegraphics[width=1.2\linewidth]{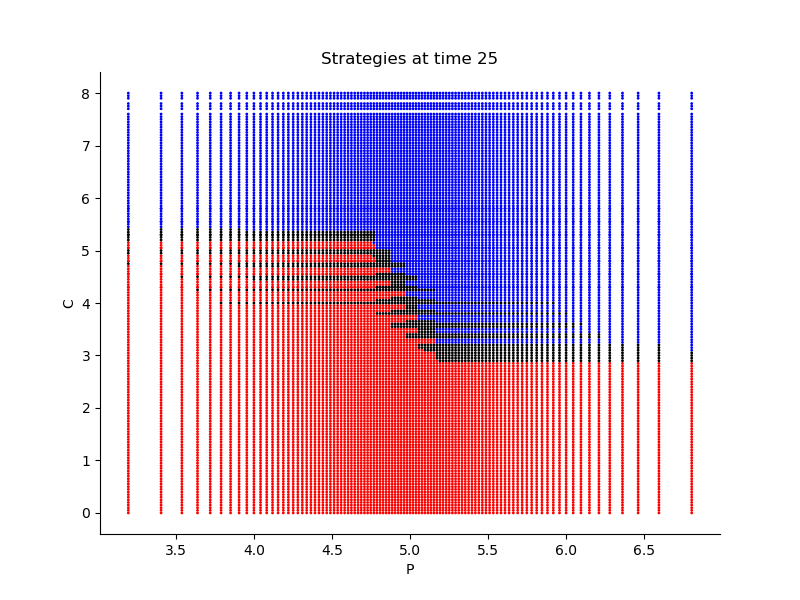}
	}
	\vspace*{-.8cm}
	\caption*{}%
\end{subfigure}
\begin{subfigure}{.5\linewidth}\centering
	{	
		\includegraphics[width=1.2\linewidth]{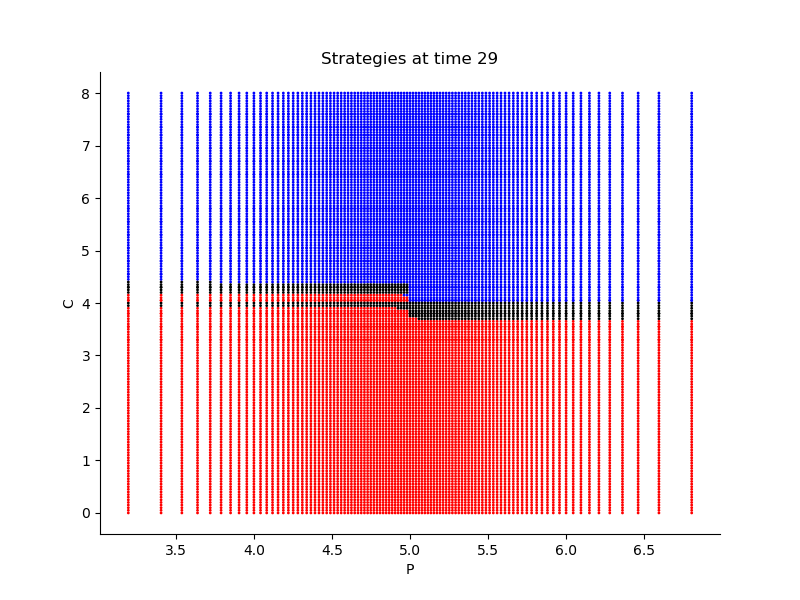}
	}
	\vspace*{-.8cm}
	\caption*{}%
\end{subfigure}
	\caption{\small{Estimated optimal decisions at times 5, 10, 15, 20, 25, 29 w.r.t. (P,C) for the  energy storage valuation problem using \textbf{Qknn}. Injection (a=-1) in red, store (a=0) in black and withdraw (a=1) in blue.}}
	\label{fig:VESRegionControl}
\end{figure}

\begin{figure}[H]
	\begin{subfigure}{.5\linewidth}\centering
		{
			\includegraphics[width=1.2\linewidth]{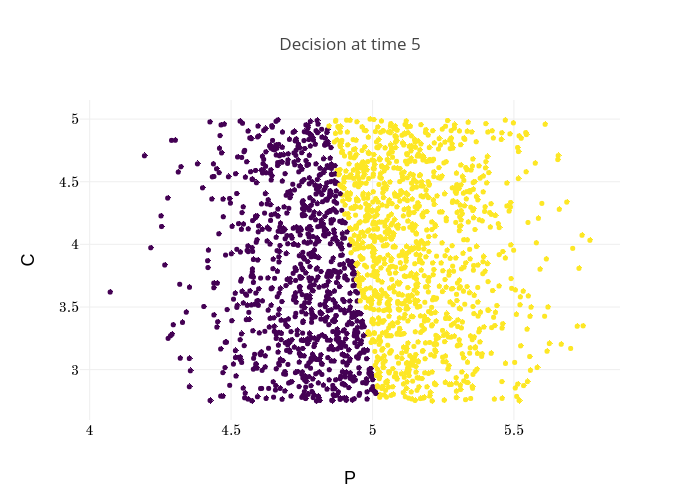}
		}
		\vspace*{-.8cm}
		\caption*{}%
	\end{subfigure}%
	\begin{subfigure}{.5\linewidth}\centering
		{	
			\includegraphics[width=1.2\linewidth]{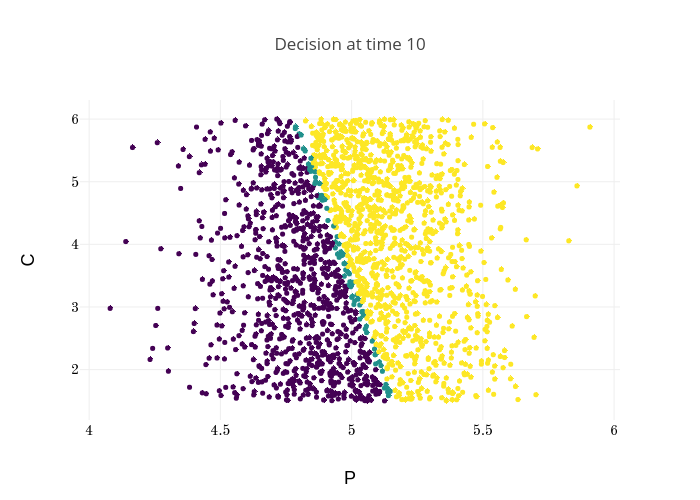}
		}
		\vspace*{-.8cm}
		\caption*{}%
	\end{subfigure}
	\begin{subfigure}{.5\linewidth}\centering
		{	
			\includegraphics[width=1.2\linewidth]{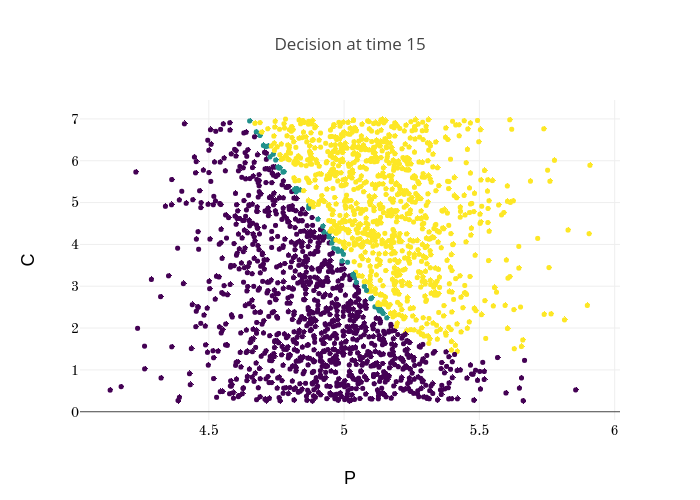}
		}
		\vspace*{-.8cm}
		\caption*{}%
	\end{subfigure}
	\begin{subfigure}{.5\linewidth}\centering
		{	
			\includegraphics[width=1.2\linewidth]{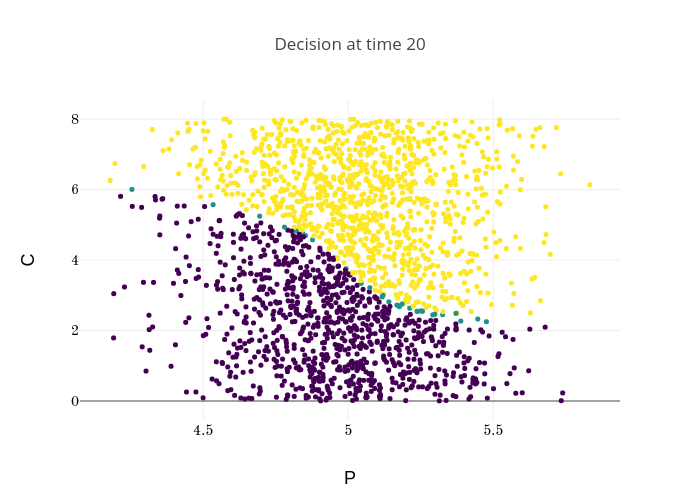}
		}
		\vspace*{-.8cm}
		\caption*{}%
	\end{subfigure}
	\begin{subfigure}{.5\linewidth}\centering
		{	
			\includegraphics[width=1.2\linewidth]{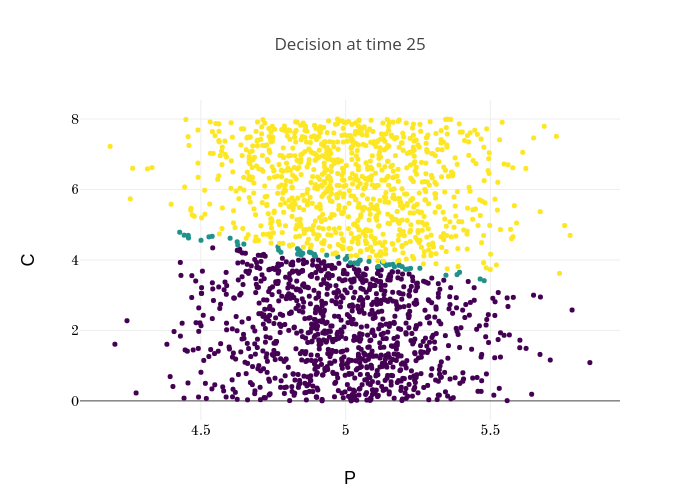}
		}
		\vspace*{-.8cm}
		\caption*{}%
	\end{subfigure}
	\begin{subfigure}{.5\linewidth}\centering
		{	
			\includegraphics[width=1.2\linewidth]{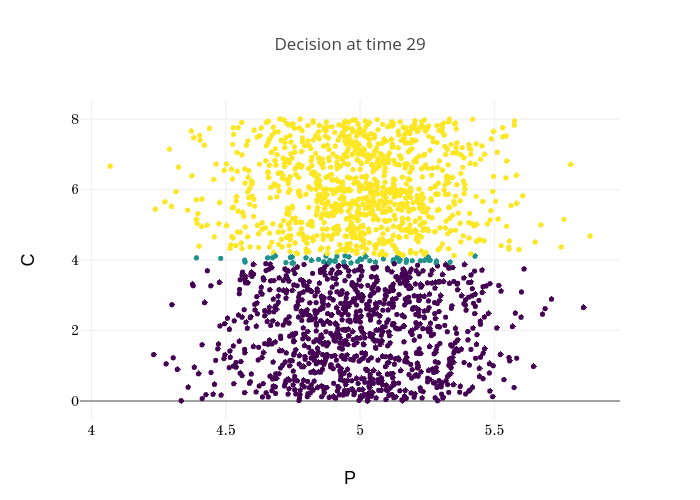}
		}
		\vspace*{-.8cm}
		\caption*{}%
	\end{subfigure}
	\caption{\small{Estimated optimal decisions at times 5, 10, 15, 20, 25, 29 w.r.t. (P,C) for the energy storage valuation problem using \textbf{ClassifPI}. Injection (a=-1) in purple, store (a=0) in blue and withdraw (a=1) in yellow.}}
	\label{fig:VESClfRegionControl}
\end{figure}

\begin{figure}[H]
	\begin{subfigure}{.5\linewidth}\centering
		{
			\includegraphics[width=1.2\linewidth]{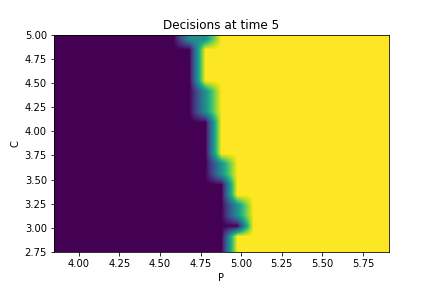}
		}
		\vspace*{-1cm}
		\caption*{}%
	\end{subfigure}%
	\begin{subfigure}{.5\linewidth}\centering
		{	
			\includegraphics[width=1.2\linewidth]{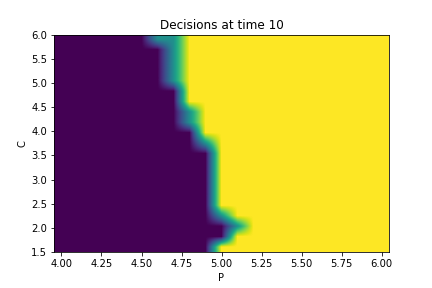}
		}
		\vspace*{-.8cm}
		\caption*{}%
	\end{subfigure}
	\begin{subfigure}{.5\linewidth}\centering
		{	
			\includegraphics[width=1.2\linewidth]{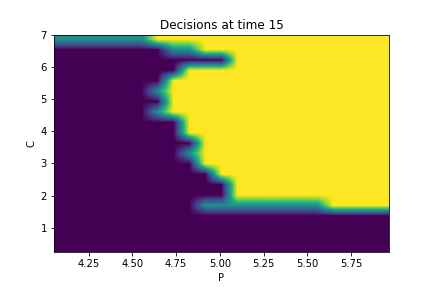}
		}
		\vspace*{-.8cm}
		\caption*{}%
	\end{subfigure}
	\begin{subfigure}{.5\linewidth}\centering
		{	
			\includegraphics[width=1.2\linewidth]{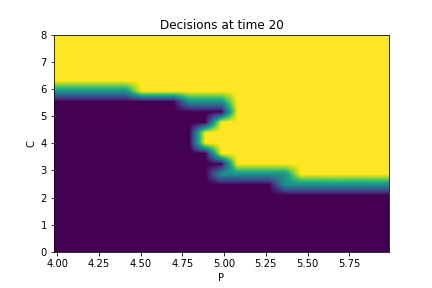}
		}
		\vspace*{-.8cm}
		\caption*{}%
	\end{subfigure}
	\begin{subfigure}{.5\linewidth}\centering
		{	
			\includegraphics[width=1.2\linewidth]{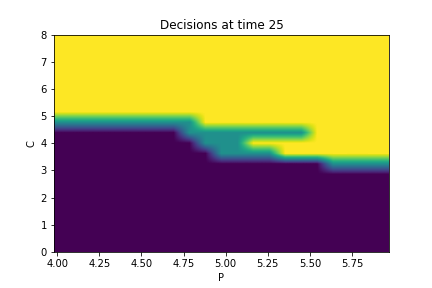}
		}
		\vspace*{-.8cm}
		\caption*{}%
	\end{subfigure}
	\begin{subfigure}{.5\linewidth}\centering
		{	
			\includegraphics[width=1.2\linewidth]{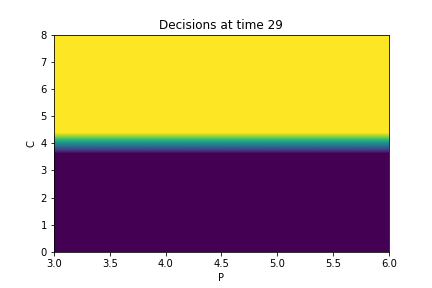}
		}
		\vspace*{-.8cm}
		\caption*{}\label{fig:VESStrat29}
	\end{subfigure}
	\caption{\small{Estimated optimal decisions at times 5, 10, 15, 20, 25, 29 w.r.t. (P,C) for the  energy storage valuation problem using \textbf{Hybrid-Now}. Injection (a=-1) in purple, store (a=0) in blue and withdraw (a=1) in yellow. Observe the instability in the decisions which come from the fact that we did not pre-train the neural networks (see Section \ref{sec:pretrain})}}
	\label{fig:VESnnRegionControl}
\end{figure}

Finally, Figures \cmt{{\color{magenta}\ref{fig:AlgINQVIRegionControl}, \ref{fig:AlgINQPIRegionControl},} } \ref{fig:VESRegionControl}, \ref{fig:VESClfRegionControl}, \ref{fig:VESnnRegionControl} provide the  optimal decisions w.r.t. $(P,C)$ at times 5, 10, 15, 20, 25, 29 estimated respectively by the Qknn, ClassifPI and Hybrid-Now algorithms. 

As expected, one can observe on each plot that the optimal strategy is to inject gas when the price is low, to sell gas when the price is high, and to make sure to have a volume of gas greater than $c_0$ in the cave when the terminal time is getting closer to minimize the terminal cost.

Let us now comment on the implementation of the algorithms: 
\begin{itemize}
\item  \textit{Qknn: }
Table \ref{t:VES} shows that once again, due to the low-dimensionality of the problem, Qknn provides the best value function estimates. The estimated optimal strategies, shown on Figure \ref{fig:VESRegionControl}, are very good estimates of the theoretical ones. The three decision regions on Figure \ref{fig:VESRegionControl} are natural and easy to interpret: basically it is optimal to sell when the price is high, and to buy when it is low. However, a closer look reveals that the waiting region (where it is optimal to do nothing) has an unusual triangular-based shape, due essentially to the discreteness of the space on which the $C$ component of the state space takes its values. We expect this shape to be very hard to reproduce with the DNN-based algorithms proposed in Section \ref{secalgo}.
\item  \textit{ClassifPI:} As shown on Figure \ref{fig:VESClfRegionControl}, the ClassifPI algorithm manages to provide accurate estimates for the optimal controls at time $n=0,\ldots,N-1$. However, the latter is not able to catch the particular triangular-based shape of the waiting region, which explains why Qknn performs better.
\item  \textit{Hybrid-Now:} As shown on Figure \ref{fig:VESnnRegionControl}, Hybrid-Now only manages to provide relatively poor estimates, compared to ClassifPI and Qknn, of the three different regions at time $n=0,\ldots,N-1$. In particular, the regions suffer from instability. 
\end{itemize}

\noindent We end this paragraph by providing some implementation details for the different algorithms we tested.

\begin{itemize}
\item  \textit{Qknn:} We used the extension of Algorithm \ref{algo:Qknn} introduced in the paragraph ``semi-linear interpolation'' of the Section 3.2.2. in \cite{BalHurPha18} and used a projection of each state on its $k$=2-nearest neighbors to get an estimate of the value function which is continuous w.r.t. the control variable at each time $n=0,\ldots,N-1$. The optimal control is computed at each point of the grids using the Brent algorithm, which is a deterministic function optimizer already implemented in Python\footnote{We could have chosen other algorithms to optimize  the $Q$-value, but, in our tests, Brent was faster than the other choices that we tried, such as GoldenSearch, and always provided accurate estimates of the optimal controls.}.
\item  \textit{Implementation details for the neural network-based algorithms:}
We use neural networks with two hidden layers, ELU activation functions\footnote{The Exponential Linear Unit (ELU) activation function is defined as $x \mapsto 
	\left\{ 
	\begin{array}{ll}
	\exp(x)-1 & \text{ if } x \le 0 \\
	x & \text{ if } x>0
	\end{array}
	\right.
	 $. }  and $20+20$ neurons \cmt{{\color{magenta} for Qknn, ClfPI and Qnn and $60+60$ neurons for  AlgINQVI and  AlgINQPI}}. The output layer contains 3 neurons with softmax activation function for the ClassifPI algorithm and no activation function for the Hybrid-Now one. \cmt{{\color{magenta} We use also  softmax activation function in  the output layer for AlgINQVI and  AlgINQPI.}}  We use a training set of size M=60,000 at each time step. Note that given the expression of the terminal cost, the ReLU activation functions (Rectified Linear Units) could have been deemed a better choice to capture the shape of the value functions, but our tests revealed that ELU activation functions provide better results. 
At time $n=0,\ldots,N-1$, we took $\mu_n=\mathcal{U}(C_{min},C_{max})$ as training measure. 

We did not use the pre-train trick discussed in Section \ref{sec:pretrain}, which explains the instability in the decisions that can be observed in Figure \ref{fig:VESnnRegionControl}.
\end{itemize}

The main conclusion of our numerical comparisons on this energy storage example is that ClassifPI, the DNN-based classification algorithm designed for stochastic control problems with discrete control space, appears to be 
more accurate than the more general Hybrid-Now. Nevertheless, ClassifPI was  not able to capture the unusual triangle-based shape of the optimal control as well as Qknn did.

\subsection{Microgrid management} 
\label{sec:sgm}

Finally, we consider a discrete-time model for power microgrid inspired by the continuous-time models developed in \cite{heyetal15} and \cite{jiapow15}; see also 
\cite{aletal18}.  
The microgrid consists of a photovoltaic (PV) power plant,  a diesel generator and a battery energy storage system (BES), hence using  a mix of  fuel and renewable energy sources. These generation units are decentralized, i.e.,  installed at a rather small scale (a few kW power), and physically close  to electricity consumers. The PV produces elec\-tricity from solar panels with a generation pattern $(P_n)_n$ depending on the weather conditions. The diesel generator has two modes: on and off. Turning it on consumes fuel, and produces an amount of power $\alpha_n$. 
The BES can store energy for later use but has limited capacity and power.  The aim of the microgrid management is to find the optimal planning that meets the power demand, denoted by 
$(D_n)_n$,  while minimizing the operational costs due to the diesel generator.  We denote by 
\beqs
R_n &=& D_n - P_n, 
\enqs
the residual demand of power: when $R_n$ $>$ $0$, one should provide power through diesel or battery, and when $R_n$ $<$ $0$, one can store the surplus power in the battery. 

The optimal control problem over a fixed horizon $N$ is formulated as follows. At any time $n = 0,\ldots,N-1$, the microgrid manager decides the power production of the diesel generator, either by 
turning it off: $\alpha_n$ $=$ $0$, or by turning it on, hence generating a power $\alpha_n$ valued in $[A_{min},A_{max}]$ with $0$ $<$ $A_{min}$ $<$ $A_{max}$ $<$ $\infty$. 
There is a fixed cost $\kappa$ $>$ $0$ associated with switching from the on/off mode to the other one off/on, and 
we denote by $M_n^\alpha$ the mode valued in $\{0 = \mbox{off}, 1= \mbox{on}\}$  of the generator right before time $n$, i.e., $M_{n+1}^\alpha$ $=$ $1_{\alpha_{n} \neq 0}$.  
 
When  the diesel generator and renewable provide a surplus of power, the excess can be stored into the battery (up to its limited capacity) for later use, and in case of power insufficiency, the battery is discharged for satisfying the power demand.   The input power process $\Ic^\alpha$ 
for charging the battery is then given by 
\beqs
\Ic_n^\alpha & = & ( \alpha_n  -  R_n)_+ \wedge (C_{max} - C_n^\alpha), 
\enqs
where $C_{max}$ is the maximum capacity of the battery with current charge $C^\alpha$, while the output power process $O^\alpha$ 
for discharging the battery is  given by 
\beqs
O_n^\alpha & = & (R_n -  \alpha_n)_+ \wedge C_n^\alpha.
\enqs
Here,  we denote $p_+$ $=$ $\max(p,0)$.  
Assuming for simplicity that the battery is fully efficient,  the capacity charge $(C_n^\alpha)_n$ of the BES, valued  in $[0,C_{max}]$, 
evolves  according to the dynamics
\beq \label{sgm::dynC}
C_{n+1}^\alpha & = &  C_n^\alpha +  \Ic_n^\alpha  - O_n^\alpha.
\enq
The imbalance process defined by 
\beqs
S_n^\alpha &=& R_n - \alpha_n + \Ic_n^\alpha - O_n^\alpha
\enqs
represents how well we are doing for satisfying electricity supply: the ideal situation occurs when $S_n^\alpha$ $=$ $0$, i.e., perfect balance between demand and generation. 
When $S_n^\alpha$ $>$ $0$, this means that demand is not satisfied, i.e., there is missing power in the microgrid, and when $S_n^\alpha$ $<$ $0$, there is an excess of electricity.  
In order to ensure that there is no missing power, we impose the following constraint on the admissible control: 
\beqs
S_n^\alpha \;  \leq \;  0, & \mbox{ i.e. } &  \alpha_n  \; \geq \; R_n   - C_n^\alpha, 
\enqs
but  penalize the excess of electricity when $S_n^\alpha$ $<$ $0$ with a proportional cost $Q^-$ $>$ $0$. We model the residual demand as a mean-reverting process:
\beqs
R_{n+1} &=&  \bar R(1-\varrho) + \varrho R_n + \eps_{n+1},
\enqs
where $(\eps_n)_n$ are i.i.d., $\bar R$ $\in$ $\R$, and $\varrho$ $<$ $1$. 
The goal of the microgrid manager is to find the optimal (admissible) decision $\alpha$ that minimizes the functional cost
\beqs
J(\alpha) &=& \E \left[ \sum_{n=0}^{N-1} \ell(\alpha_n) +  \kappa 1_{\{M_n^\alpha \neq M_{n+1}^\alpha\}} + Q^- (S_n^\alpha)_-   \right], 
\enqs
where $\ell(.)$ is the cost function for fuel consumption: $\ell(0)$ $=$ $0$, and e.g. $\ell(a)$ $=$ $K a^\gamma$, with $K$ $>$ $0$, $\gamma$ $>$ $0$. 
This stochastic control problem fits into the $3$-dimensional framework of Section \ref{secintro} (see also Remark \ref{remcons}) with control $\alpha$  valued in $\A$ $=$ $\{0\}\times [A_{min},A_{max}]$, $X^\alpha$ $=$ $(C^\alpha,M^\alpha,R)$, noise 
$\eps_{n+1}$, starting from an initial value  $(C_0^\alpha,M_0^\alpha,R_0)$ $=$ $(c_0,0,r_0)$  on the state space  
$[0,C_{max}]\times\{0,1\}\times\R$,  with dynamics function    
\beqs
F(x,a,e) & = & 
\left(
\begin{array}{c}
F^1(x,a) := c +  (a  -  r)_+ \wedge (C_{max} - c)  -  (r - a)_+  \wedge c    \\
1_{a \neq 0}  \\
\bar R(1-\varrho) + \varrho r  + e 
\end{array}
\right),
\enqs
for $x$ $=$ $(c,m,r)$ $\in$ $[0,C_{max}]\times\{0,1\}\times\R$, $a$ $\in$ $\{0\}\times [A_{min},A_{max}]$, 
$e$  $\in$ $\R$,  running cost function
\beqs
f(x,a) & = & \ell(a)  +  \kappa 1_{ m = 1_{a=0}} + Q^- S(x,a)_-,    \\
S(x,a) &=&  r  - a  +  (a  - r)_+ \wedge (C_{max} - c)  -  (r- a)_+  \wedge c,  
\enqs
zero terminal cost $g$ $=$ $0$, and control  constraint
\beqs
\A_n(x) & = & \Big\{ a  \in \{0\}\times [A_{min},A_{max}]:  S(x,a) \leq 0  \Big\} \\
&=&   \Big\{ a  \in \{0\}\times [A_{min},A_{max}]:  r  - c    \leq  a   \Big\}. 
\enqs

\begin{Remark}
{\rm	The state/space constraint is managed in our NN-based algorithm by introducing a penalty function into the running cost (see Remark \ref{remcons}): $f(x,a)$ $\leftarrow$ $f(x,a)+L(x,a)$ 
	\beqs
	L(x,a) &=& Q^+ \Big(  r  - c    - a  \Big)_+
	\enqs  
	with large $Q^+$ taken much larger than $Q^-$. Doing so, the NN-based estimate of the optimal control learns not to take any forbidden decision.
	}
\ep
\end{Remark}

The control space $\{0\} \cup [A_\mathrm{min}, A_\mathrm{max}]$ is a mix between a discrete space and a continuous space, which is challenging for algorithms with neural networks. We actually use 
a mixture of classification and standard DNN  for the control: $(p_0(x;\theta),\pi(x;\beta))$ valued in $[0,1]\times [A_{min},A_{max}]$, where  
$p_0(x;\theta)$ is the probability of turning off in state $x$, and $\pi(x;\beta)$ is the amount of power when turning on with probability $1-p_0(x;\theta)$.   In other words, 
\beqs
X_{n+1} 
&=&
\left\{
 \begin{array}{ll}
F(X_n,0,\eps_{n+1}) & \mbox{ with probability } p_0(X_n;\theta_n)  \\
F(X_n,\pi(X_n;\beta_n),\eps_{n+1}) & \mbox{ with probability } 1- p_0(X_n;\theta_n) 
\end{array}
\right.
\enqs
The pseudo-code of this approach, specifically designed for this problem, is written in Algorithm \ref{algo:ClassifHybrid}, and we henceforth refer to it as ClassifHybrid. Note in particular that it is an Hybrid version of ClassifPI. 

\vspace{3mm}
\begin{algorithm}[H]
	\caption{ClassifHybrid}
	\label{algo:ClassifHybrid}
	\textbf{Input:} the training distributions $(\mu_n)_{n=0}^{N-1}$\;
	\textbf{Output:} \\
	-- estimate of the optimal strategy $(\hat a_n)_{n=0}^{N-1}$\;
	-- estimate of the value function $(\hat V_n)_{n=0}^{N-1}$\;
	Set $\hat V_N$ $=$ $g$\;
	\For{$n$ $=$ $N-1,\ldots,0$}{
Compute
\begin{align*}
(\hat\beta_n^0,\hat\beta_n^1) &\in 
\argmax_{\beta^0,\beta^1} \E \Bigg[ p_0(X_n;\beta^0) \left[ f(X_n,0)  + \hat V_{n+1} \left(f(\hat X_{n+1}^{0}  \right) \right]\\
&  + (1- p_0(X_n;\beta^0) )  \left[ f(X_n,\pi(X_n;\beta^1))  + 
\hat V_{n+1} \left(\hat X_{n+1}^{1,\beta^1} \right) \right] \Bigg], 
\end{align*}
where $X_n$ $\leadsto$ $\mu_n$, $\hat X_{n+1}^{0}$ $=$ $F(X_n,0,\eps_{n+1})$,  and
$\hat X_{n+1}^{1,\beta^1}$ $=$ $F(X_n,\pi(X_n;\beta^1),\eps_{n+1})$\;
		Compute
\begin{align*}
\hat\theta_n  &\in 
\argmin_{\theta} \E \Bigg[ p_0\big(X_n;\hat \beta_n^0\big) \left[ f(X_n,0)  + \hat V_{n+1} \left(f(\hat X_{n+1}^{0}  \right) -\Phi(.;\theta) \right]^2 \\
&  + \big(1- p_0\big(X_n;\hat \beta_n^0\big) \big)  \left[ f(X_n,\pi(X_n;\beta_n^1))  + 
\hat V_{n+1} \left(\hat X_{n+1}^{1,\hat \beta_n^1} \right) -\Phi(.;\theta) \right]^2 \Bigg];
\end{align*}
		Set $\hat V_n = \Phi(.;\hat\theta_n) $; \Comment{$\hat V_n$ is the estimate of the value function at time $n$}
	}
\end{algorithm}

\vspace{2mm}

\paragraph{Test}
We set the parameters to the following values to compare Qknn and ClassifHybrid:
\begin{equation*}
\begin{array}{rclrclrclrcl}
N&=&30 \text{ or } 200,&\qquad   \bar R &=&0.1,&\qquad   \varrho&=&0.9,&\qquad   \sigma&=&0.2,\\ C_\mathrm{min}&=&0,& C_\mathrm{max}&=&1 \text{ or } 4, & C_0&=&0,&K&=&2,\\ \gamma&=&2, &\kappa&=&0.2,& Q^-&=&10,& R_0&=&0.1,\\ A_\mathrm{min}&=&0.05,& A_\mathrm{max}&=&10 & Q^+&=&1000.
\end{array}
\end{equation*}

\begin{figure}[H]
	\centering
	\begin{subfigure}[b]{1\textwidth}
		\includegraphics[width=1.1\linewidth]{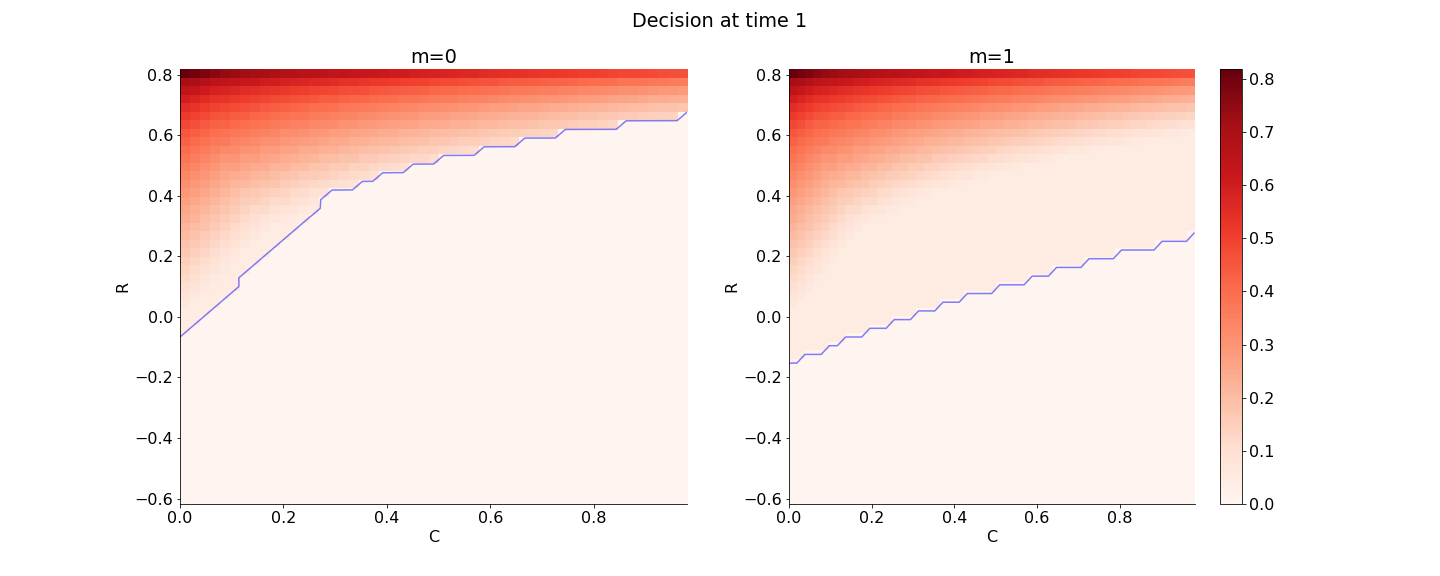}
		\caption*{}
		\vspace*{-1cm}
	\end{subfigure}
	
	\begin{subfigure}[b]{1\textwidth}
		\includegraphics[width=1.1\linewidth]{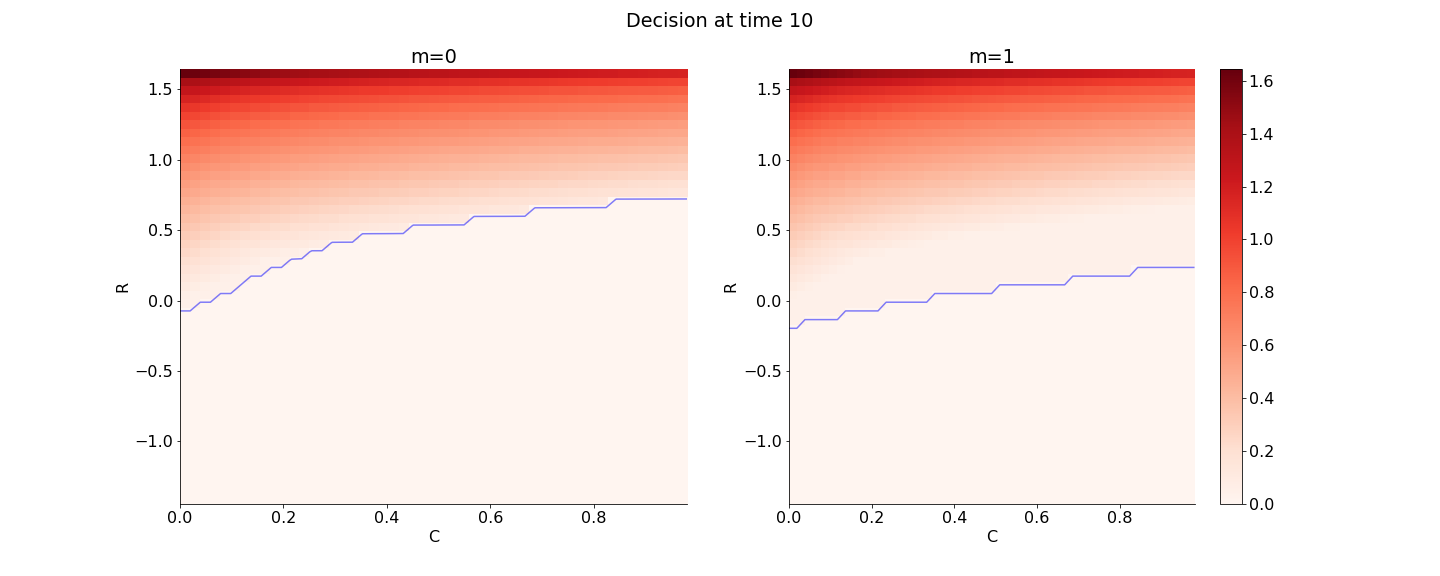}
		\caption*{}
		\vspace*{-1cm}
	\end{subfigure}
	
	\begin{subfigure}[b]{1\textwidth}
		\includegraphics[width=1.1\linewidth]{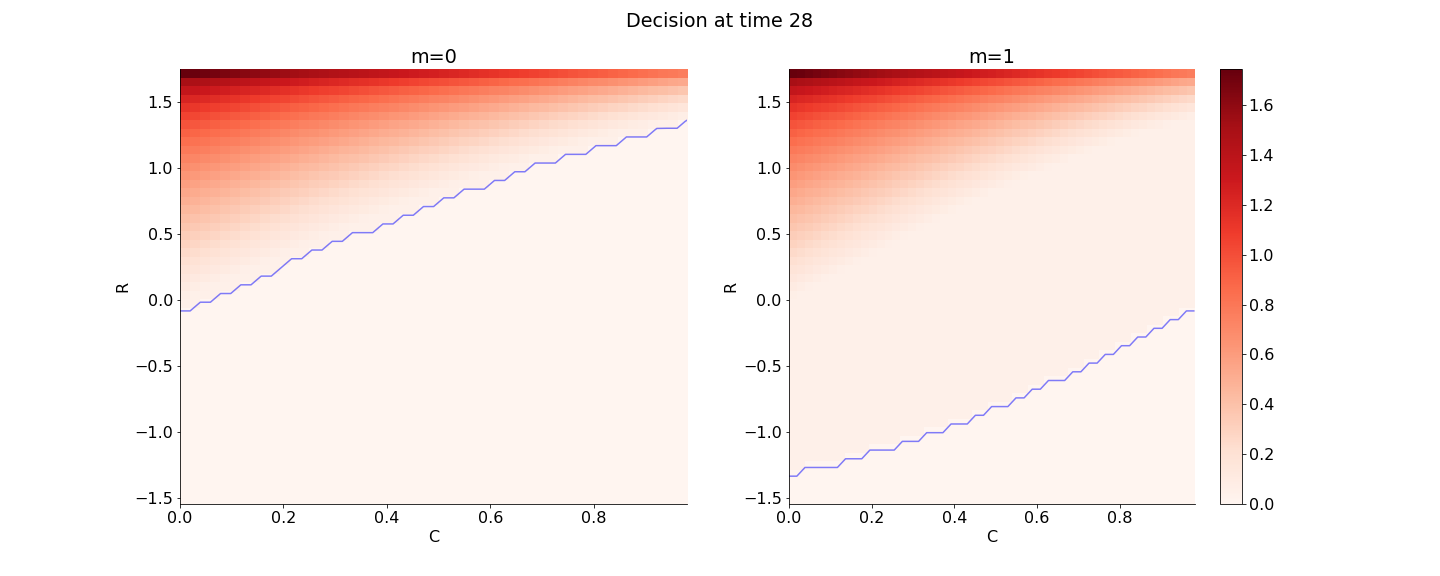}
		\caption*{}
		\vspace*{-1cm}
	\end{subfigure}
	\caption{\small{Estimated optimal decisions at time 1, 10 and 28, using Qknn, with $N=30$ time steps. The region under the blue line is the one where it is optimal to turn the generator off if $m$=1 (i.e. the generator was on at time $n$-1), or keep it off if $m=0$ (i.e. the generator was off at time $n$-1).}}
	\label{fig:Qknn_decisions_SGM}
\end{figure}

\begin{figure}[H]
	\centering
	\begin{subfigure}[b]{1\textwidth}
		\includegraphics[width=1.1\linewidth]{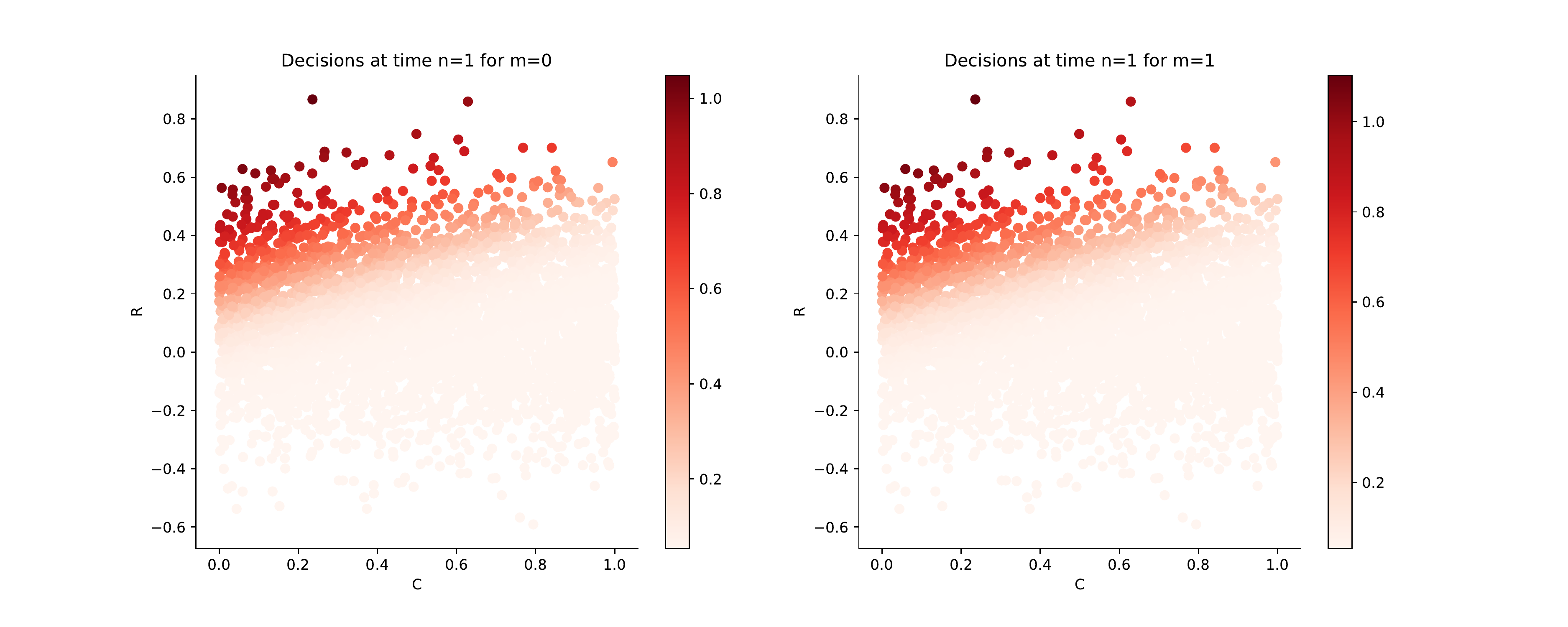}
		\caption*{}
		\vspace*{-1cm}
	\end{subfigure}
	
	\begin{subfigure}[b]{1\textwidth}
		\includegraphics[width=1.1\linewidth]{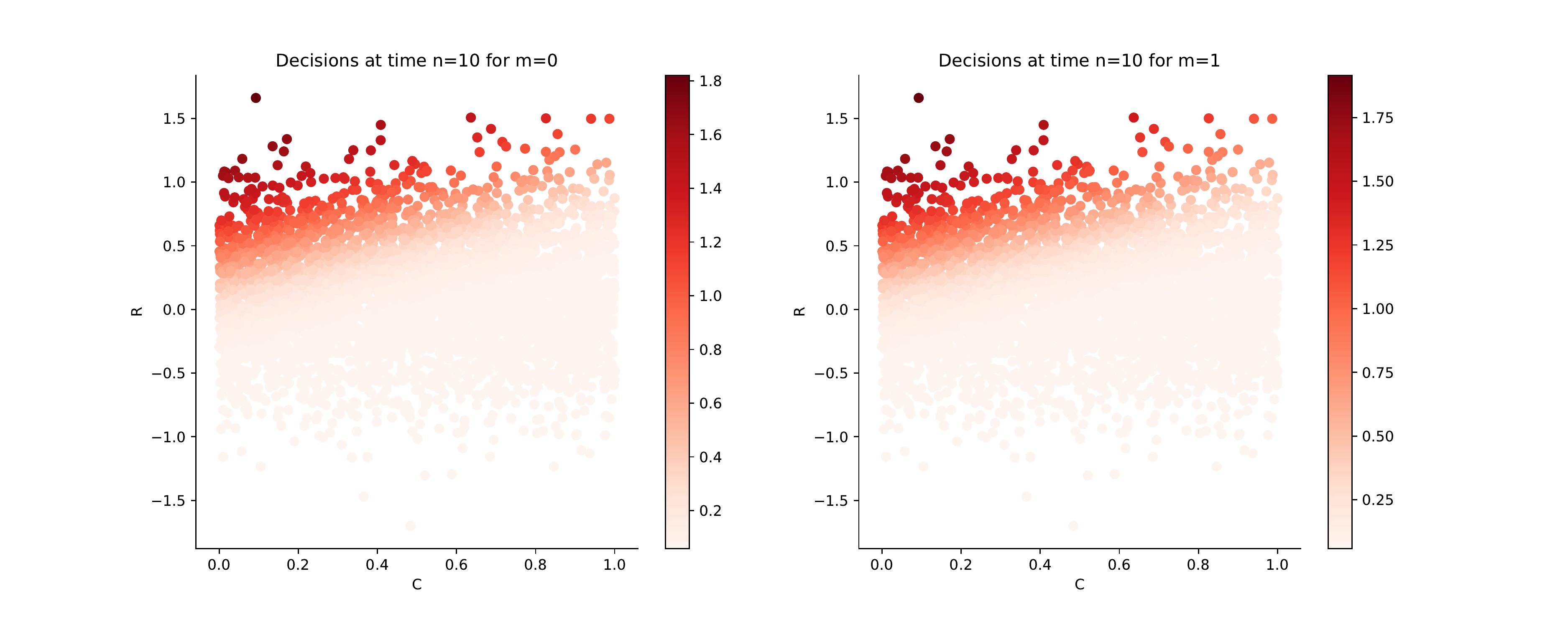}
		\caption*{}
		\vspace*{-1cm}
	\end{subfigure}
	
	\begin{subfigure}[b]{1\textwidth}
		\includegraphics[width=1.1\linewidth]{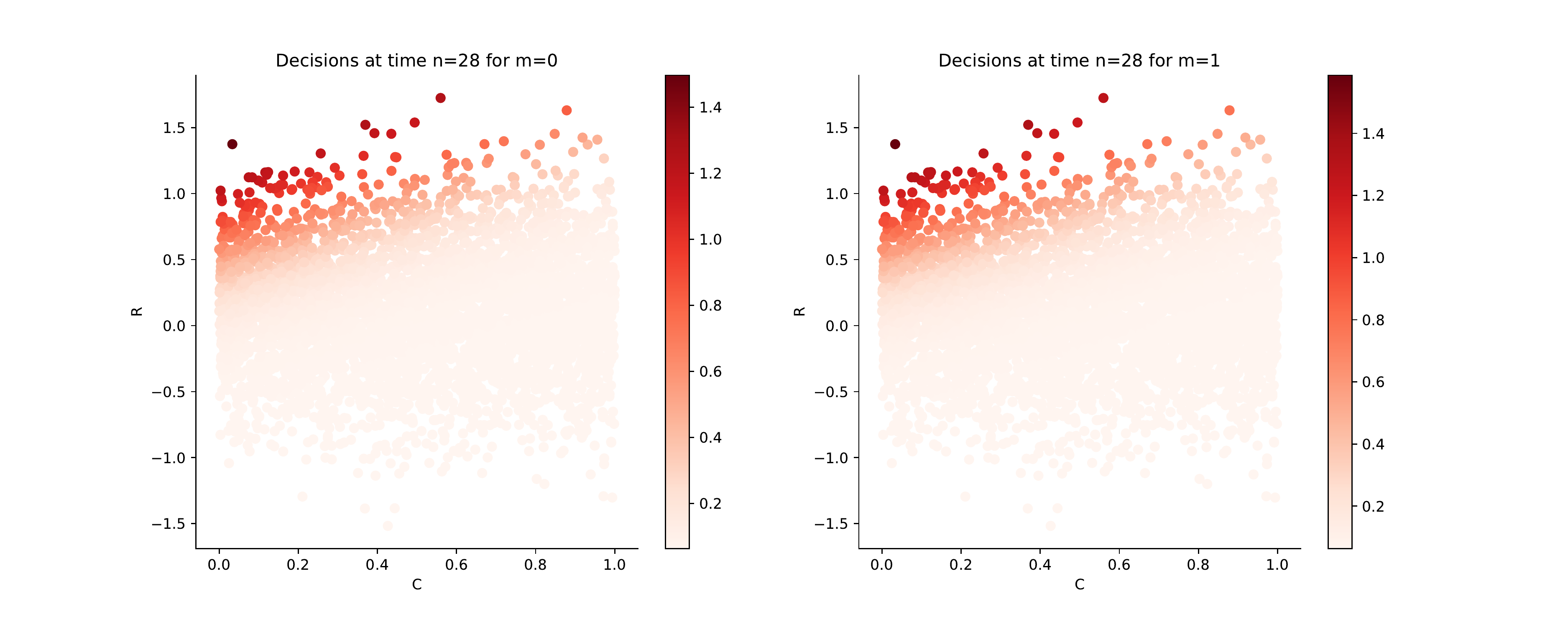}
		\caption*{}
		\vspace*{-1cm}
	\end{subfigure}
	\caption{\small{Estimated optimal decisions at time 1, 10 and 28, using ClassifHybrid, with $N=30$ time steps. }}
	\label{fig:classifVI_decisions_SGM}
\end{figure}

\paragraph{Results}
Figure \ref{fig:Qknn_decisions_SGM} shows the Qknn-estimated optimal decisions to take at times $n$ $=$ $1, 10, 28$ in the cases where $m$ $=$ $M_n$ $=$ $0$ and  $m$ $=$ $M_n$ $=$ $1$. 
If the generator is off at time $n$, i.e. $m$ $=$ $0$, the blue curve separates the region where it is optimal to keep it off and the one where it is optimal to generate power. If the generator is on at time $n$, i.e. $m=1$, the blue curve separates the region where it is optimal to turn it off and the one where it is optimal to generate power. A colorscale is available on the right to inform how much power it is optimal to generate in both cases. Observe that the optimal decisions are quite intuitive: for example,  if the demand is high and the battery is empty, then it is optimal to generate a lot of energy. Moreover,  it is optimal to turn the generator off if the demand is negative or if the battery is charged enough to meet the demand. \\
We plot in Figure \ref{fig:classifVI_decisions_SGM}  the estimated optimal decisions at times $n$ $=$ $1, 10, 28$, using the Hybrid-Now algorithm, with $N=30$ time steps. See that the decisions are similar to the ones given using Qknn.\\
Note that the plots in Figure \ref{fig:Qknn_decisions_SGM} and \ref{fig:classifVI_decisions_SGM} look much better than the ones obtained in \cite{aletal18} in which algorithms based on regress-now or regress-later are used (see in particular Figure 4 in \cite{aletal18}); hence Qknn and ClassifHybrid seem more stable than the algorithms proposed in \cite{aletal18}. 

We report in Table \ref{t:SmartGridManagement2} the result for the estimates of the value function with $N$=30 time steps, obtained by running 10 times a forward Monte Carlo with 10,000 simulations using the optimal strategy estimated using Qknn and ClassifHybrid algorithms. Observe that Hybrid-Now performs better than Qknn. However, Qknn run in less than a minute whereas Hybrid-Now needed seven minutes to run.\\
We also report in Table \ref{t:SmartGridManagement} the value function estimates with $N$=200 time steps, obtained by running 20 times a forward Monte Carlo with 10,000 simulations using the Qknn-estimated optimal strategy. 

\begin{table}[H]
	\centering
	\caption{\small{Estimates of the value function at time 0 and state $(C_0=0,M_0=0,R_0=0.1)$, for $N=30$ and $C_{max}=1$, using Qknn and ClassifHybrid algorithms. Note that ClassifHybrid achieved better results than Qknn on this problem.}}
	\label{t:SmartGridManagement2}
	\vspace{-3mm}
	\begin{tabular}{c|c|c}
	&	Mean   &   std    \\ \hline
	ClassifHybrid &	33.34 & 0.31\\
	Qknn &	35.37 & 0.34
	\end{tabular}
\end{table}

\begin{table}[H]
	\centering
	\caption{\small{Qknn-estimates of the value function at time 0 and state $(C_0=0,M_0=0,R_0=0.1)$, for $N=200$.}}
	\label{t:SmartGridManagement}
	\vspace{-3mm}
	\begin{tabular}{c|c}
 Mean   &   Standard Deviation  \\ \hline
231.8 & 1.2
	\end{tabular}
\end{table}

\begin{figure}[H]
	\centering
	\begin{subfigure}[b]{1\textwidth}
		\includegraphics[width=1.1\linewidth]{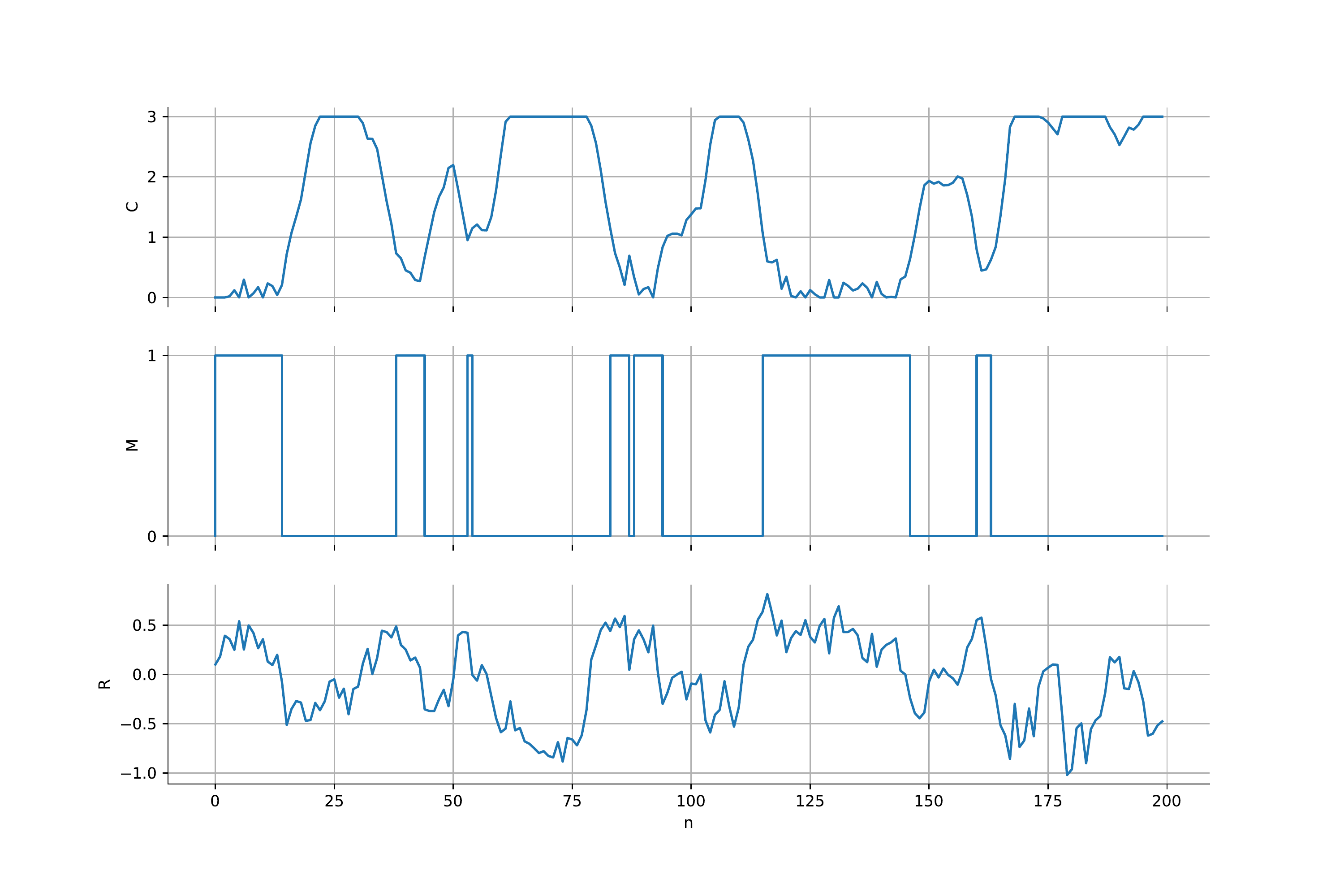}
		\caption*{}
		\vspace*{-1cm}
	\end{subfigure}	
	\begin{subfigure}[b]{1\textwidth}
		\includegraphics[width=1.1\linewidth]{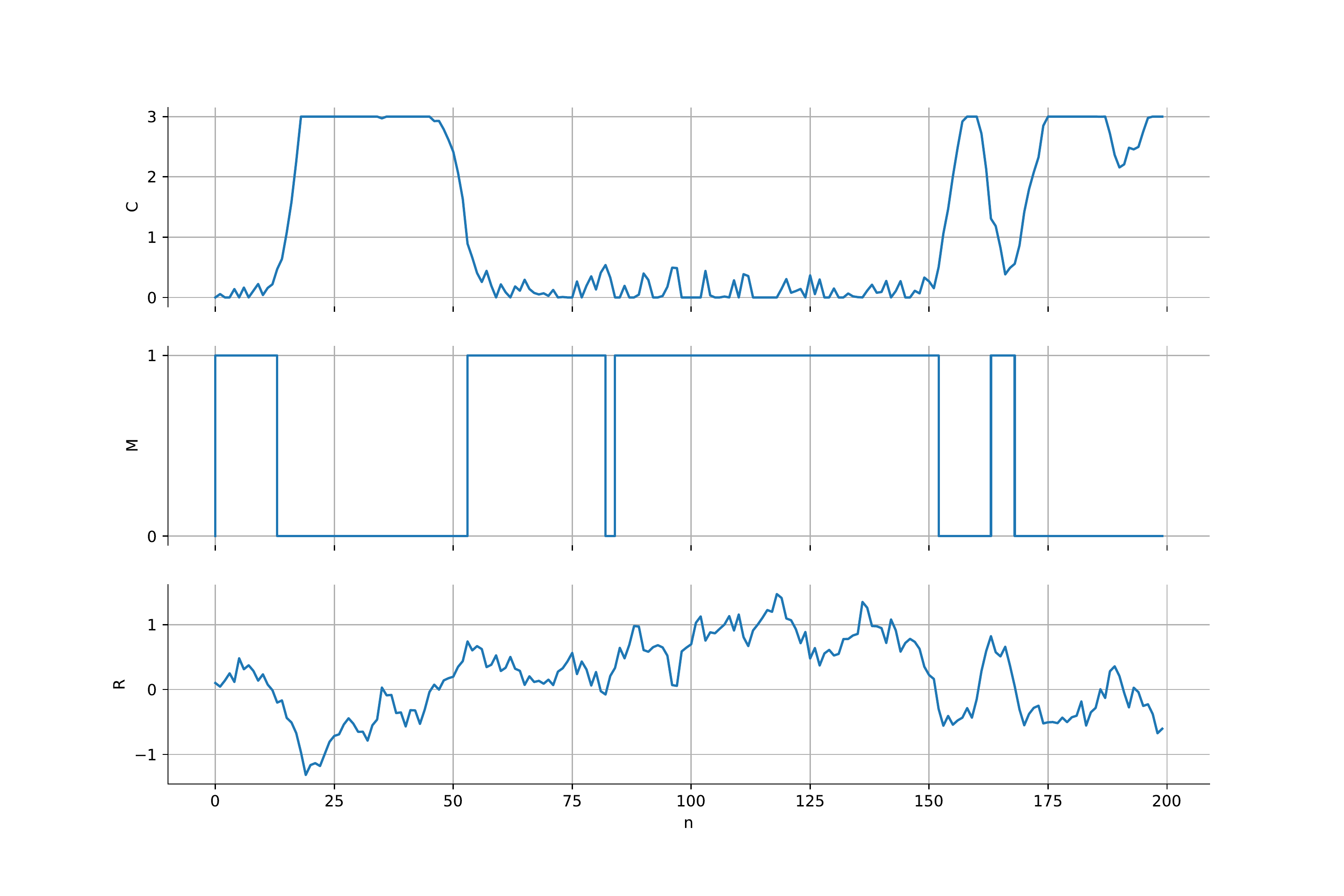}
		\caption*{}
		\vspace*{-1cm}
	\end{subfigure}
	\caption{\small{Two simulations of $(C,M,R)$ optimally controlled using Qknn, with $N=200$ and $C_\mathrm{max}=4$.}}
	\label{fig:Qknn_simuls_SGM}
\end{figure}

Figure \ref{fig:Qknn_simuls_SGM} shows two simulations of $(C,M,R)$ controlled using the Qknn-estimated optimal strategy, where $N=200$ has been chosen. Observe in particular the natural behavior of the Qknn-decisions which consists in turning the generator on when the demand cannot be met by the battery, and turn it off when the demand is negative or when the battery is charged enough to meet the demand. Note that the plots are similar to the ones plotted in Figure 9 of \cite{aletal18}.

\vspace{2mm}
\noindent \textbf{Comments on Qknn:} Note that there is no need to use a penalization method with the Qknn-algorithm to constrain the control to stay in $\mathbb{A}_n(x)$, where $x$ is the state at time $n$, since, for all state $x$, we can simply search for the optimal control associated in $\mathbb{A}_n(x)$, using e.g. the Brent algorithm. For $n=0,\ldots,N-1$, we took the training set as follows: $\Gamma_n:= \Gamma_C \times \{0,1\} \times \Gamma_R^n$; where $\Gamma_C:= \{C_{min} + \frac{i}{50} (C_{max}-C_{min} ), i=0,\ldots, 50\} $, $\Gamma_R^n:=\rho^n R_0+ \sigma \frac{1-\rho^n}{1-\rho}   \Gamma_1$ and where $\Gamma_1$ is the optimal grid for the quantization of $\mathcal{N}(0,1)$, available in \url{http://www.quantize.maths-fi.com}, with 51 points. This choice of training points for the $C$ component corresponds to the exploration procedure discussed in Remark \ref{rk:trainingMeasure}, whereas we chose the best grid with 51 points for the (uncontrolled) $R$ component.\\
\textbf{Comments on ClassifHybrid:} We took 100 mini-batches of size 300 and took 100 epochs to run the algorithm. We chose the following training distribution at time $n$: $\mu_n= \mathcal{U}(C_{min},C_{max}) \times \mathcal{U}(\{0,1\}) \times \P_{R_n}$, where $\P_{R_n}$ is the law of the (uncontrolled) residual demand at time $t_n$. Note that such a choice of training distribution means that we want to explore all the available states for the controlled components of the controlled process $(C,M,R)$ in order to learn the optimal strategy globally.\\

The microgrid management problem is very challenging for our algorithms because the control space $\{0\} \cup [a_\mathrm{min}, a_\mathrm{max}]$ is a mix of discrete and continuous space, moreover the choice of the optimal control is subject to constraints. We designed ClassifHybrid, an Hybrid version of ClassifPI, to solve this problem. ClassifHybrid provided very good estimates and actually managed to perform better than Qknn. 

\section{Discussion and conclusion} \label{secconclu}

Our proposed algorithms are well-designed and provide accurate estimates of optimal control and value function associated with various high-dimensional control problems. 
Also, when tested on low-dimensional problems, they performed as well as the Monte Carlo-based or quantization-based methods,  which have shown their efficiency in low dimension, see e.g. \cite{BalHurPha18} and \cite{aletal18}. 
 
The presented algorithms suffer from a rather high time-consuming cost due to the expensive training of $2(N-1)$ neural networks to learn the value functions and optimal controls at times $n$ $=$ $0,\ldots,N-1$. 
However, the agent can easily alleviate the computation time. A first trick consists in reducing the number of neural networks by partially or totally ignoring the dynamic programming principle (DPP), as it has been done e.g. in \cite{Ehanjen17}. The use of one unique Recurrent Neural Networks (RNN)  (in the case where the DPP is totally ignored) or a few of them (in the partial-ignored case) can also be considered to learn the optimal controls, either all at the same time (first case), or group by group in a backward way (second case). We refer to \cite{CMW18} for algorithms in this spirit.
Another trick consists in learning faster the value functions and optimal controls at times $n$ $=$ $0,\ldots,N-1$ by pre-training the neural networks. The way to proceed in that direction is to initialize at time $n$ the weights and bias of the value function estimator $\hat{V}_n$ to the ones of $\hat{V}_{n+1}$. We then rely on the continuity of the value function w.r.t. the time $n$ to expect that the weights will not change much from time $n$ to $n+1$, hence trainable very quickly by reducing the learning rate of the Adam algorithm for the gradient descent, and using an early-stop procedure as implemented in Keras\footnote{See EarlyStopping callback in Keras}. Another benefit from the pre-training task is to get the stability of the estimates w.r.t. time, which is also a pleasant feature.

\small

\printbibliography

\end{document}